\RequirePackage{ifpdf}
\ifpdf 
\documentclass[pdftex]{sigma}
\else
\documentclass{sigma}
\fi

\begin{document}


\renewcommand{\thefootnote}{$\star$}

\renewcommand{\PaperNumber}{062}

\FirstPageHeading

\ShortArticleName{Isoparametric and Dupin Hypersurfaces}

\ArticleName{Isoparametric and Dupin Hypersurfaces\footnote{This paper is a
contribution to the Special Issue ``\'Elie Cartan and Dif\/ferential Geometry''. The
full collection is available at
\href{http://www.emis.de/journals/SIGMA/Cartan.html}{http://www.emis.de/journals/SIGMA/Cartan.html}}}

\Author{Thomas E. CECIL}

\AuthorNameForHeading{T.E. Cecil}

\Address{Department of Mathematics and Computer Science, College of the Holy Cross,\\
Worcester, MA 01610, USA}
\Email{\href{mailto:cecil@mathcs.holycross.edu}{cecil@mathcs.holycross.edu}}
\URLaddress{\url{http://mathcs.holycross.edu/faculty/cecil.html}}

\ArticleDates{Received June 24, 2008, in f\/inal form August 28,
2008; Published online September 08, 2008}

\Abstract{A hypersurface $M^{n-1}$ in a real space-form ${\bf R}^n$, $S^n$ or $H^n$ is isoparametric if it has constant principal curvatures.  For ${\bf R}^n$ and $H^n$, the classif\/ication of isoparametric hypersurfaces is complete and
relatively simple, but as \'{E}lie Cartan 
showed in a series of four papers in 1938--1940,
the subject is much deeper and more complex for hypersurfaces in the sphere $S^n$.
A hypersurface $M^{n-1}$ in a real space-form
is proper Dupin if the number $g$ of distinct principal
curvatures is constant on $M^{n-1}$, and each principal curvature function is constant along each leaf of its corresponding
principal foliation. This is an important genera\-lization of the isoparametric property that has its roots in
nineteenth century dif\/ferential geometry and has been studied ef\/fectively in the context of Lie sphere geometry.
This paper is a survey of the known results in these f\/ields with emphasis on results
that have been obtained in more recent years and discussion of important open problems in the f\/ield.}


\Classification{53C40; 53C42; 53B25}

\section{Introduction}
\label{intro}
A hypersurface $M^{n-1}$ immersed in Euclidean space ${\bf R}^n$, the sphere $S^n$ or hyperbolic space $H^n$ is
called {\em isoparametric} if it has constant principal curvatures.
An isoparametric hypersurface in ${\bf R}^n$ can have at most
two distinct principal curvatures, and it must be an open
subset of a hyperplane, hypersphere or a spherical cylinder
$S^k \times {\mathbf R}^{n-k-1}$. This was f\/irst proven for $n=3$ by Somigliana \cite{Som} in 1919
(see also B.~Segre \cite{Seg1} and Levi-Civita \cite{Lev}) and for arbitrary~$n$ by B.~Segre \cite{Seg} in~1938.
A similar result holds in $H^n$.
However, as \'{E}lie Cartan \cite{Car1,Car2,Car3,Car4} showed in a series of four papers published in the period
1938--1940, the theory of isoparametric hypersurfaces in the sphere $S^n$ is much more beautiful and complicated.

Cartan produced examples of isoparametric hypersurfaces in spheres
with $g = 1,2,3$ or 4 distinct principal
curvatures, and he classif\/ied those with $g \leq 3$. Approximately thirty years later, M\"{u}nzner
\cite{Mu,Mu2} wrote two papers
which greatly extended Cartan's work, proving that all isoparametric hypersurfaces are algebraic and that the number
$g$ of distinct principal curvatures must be  1, 2, 3, 4  or 6.  Since Cartan had classif\/ied isoparametric hypersurfaces
with $g \leq 3$, the classif\/ication of those with $g=4$ or 6 quickly became the goal of researchers in the f\/ield
after M\"{u}nzner's work.
This classif\/ication problem has proven to be interesting and dif\/f\/icult, and it was listed as Problem 34 on
Yau's \cite{Yau} list
of important open problems in geometry in 1990.  The problem remains open in both cases $g=4$ and 6 at
this time, although much progress has been made.

A hypersurface $M^{n-1}$ in one of the real space-forms ${\bf R}^n$, $S^n$ or $H^n$
is {\em proper Dupin} if the number $g$ of distinct principal
curvatures is constant on $M^{n-1}$, and each principal curvature function is constant along each leaf of its corresponding
principal foliation.  This is clearly a generalization of the isoparametric condition, and it has its roots
in nineteenth century dif\/ferential geometry.  Aside from the isoparametric surfaces in ${\bf R}^3$, the f\/irst
examples of Dupin hypersurfaces are the cyclides of Dupin in ${\bf R}^3$.  These are all surfaces that can be obtained
from a torus of revolution, a circular cylinder or a circular cone by inversion in a 2-sphere in ${\bf R}^3$.
The cyclides have several other characterizations that will be discussed in this paper (see
also \cite[pp.~151--166]{CR}, \cite[pp.~148--159]{Cec1}, \cite{P3}).

The proper Dupin property is preserved under M\"{o}bius (conformal) transformations, and an important class of
compact proper Dupin hypersurfaces in ${\bf R}^n$ consists of those hypersurfaces that are obtained
from isoparametric hypersurfaces
in a sphere $S^n$ via stereographic projection from $S^n - \{P\}$ to ${\bf R}^n$, where $P$ is a point in $S^n$.
Indeed several of the major classif\/ication results for compact
proper Dupin hypersurfaces involve such hypersurfaces.

An important step in the theory of proper Dupin hypersurfaces was the work of Pinkall~\cite{P1,P6,P2,P3}
which situated the study of Dupin hypersurfaces in the setting of Lie sphere geometry.  Among other things,
Pinkall proved that the
proper Dupin property is invariant under the group of Lie sphere transformations, which contains the
group of M\"{o}bius transformations as a subgroup.  The theory of proper Dupin hypersurfaces has both local and global
aspects to it, and many natural problems remain open.

In this paper, we survey the known results for isoparametric and proper Dupin hypersurfaces with emphasis on results
that have been obtained in more recent years, and discuss important open problems in the f\/ield.  The reader is also
referred to the excellent survey article by Thorbergsson~\cite{Th6} published in the year 2000.

\section{Isoparametric hypersurfaces}
\label{sec:1}
Let $\tilde{M}^n(c)$ be a simply connected, complete Riemannian manifold of dimension $n$ with constant sectional
curvature $c$, that is, a real space-form.  For $c= 0,1,-1$, respectively, $\tilde{M}^n(c)$ is Euclidean space
${\bf R}^n$, the unit sphere $S^n \subset {\bf R}^{n+1}$, or hyperbolic space $H^n$.  According to the original def\/inition,
a 1-parameter family $M_t$ of hypersurfaces in $\tilde{M}^n(c)$ is called an {\em isoparametric family} if each $M_t$
is equal to a level set $V^{-1}(t)$ for some non-constant smooth real-valued function $V$ def\/ined on a connected open subset
of $\tilde{M}^n(c)$ such that the gradient and Laplacian of $V$ satisfy
\begin{gather}
\label{isoparametric}
|{\rm grad}\,V|^2 = T \circ V, \qquad \triangle V = S \circ V,
\end{gather}
for some smooth functions $S$ and $T$.  Thus, the two classical Beltrami dif\/ferential parameters are both
functions of $V$ itself, which leads to the name ``isoparametric'' for such a family of hypersurfaces.
Such a function $V$ is called an {\em isoparametric function}. (See Thorbergsson \cite[pp.~965--967]{Th6}
and Q.-M.~Wang \cite{Wang3, Wang1} for more discussion of isoparametric functions.)

Let $f:M^{n-1} \rightarrow \tilde{M}^n(c)$ be an oriented hypersurface with f\/ield of unit normal vectors $\xi$.
The {\em parallel hypersurface} to $f(M^{n-1})$ at signed distance $t\in {\bf R}$ is the map
$f_t:M^{n-1} \rightarrow \tilde{M}^n(c)$ such that for each $x \in M^{n-1}$, the point $f_t(x)$ is obtained
by traveling a signed distance $t$ along the geodesic in $\tilde{M}^n(c)$ with initial point $f(x)$ and initial
tangent vector $\xi(x)$.  For $\tilde{M}^n(c) = {\bf R}^n$, the formula for $f_t$ is
\begin{gather}
\label{parallel-E}
f_t(x) = f(x) + t \xi(x),
\end{gather}
and for $\tilde{M}^n(c) = S^n$, the formula for $f_t$ is{\samepage
\begin{gather}
\label{parallel-S}
f_t(x) = \cos t \; f(x) + \sin t \; \xi(x).
\end{gather}
There is a similar formula in hyperbolic space (see, for example, \cite{CecToh}).}

Locally, for suf\/f\/iciently small values of $t$, the map $f_t$ is also an immersed hypersurface.  However, the map
$f_t$ may develop singularities at the focal points of the original hypersurface $f(M^{n-1})$.  Specif\/ically, a point
$p=f_t(x)$ is called a {\em focal point of multiplicity $m > 0$ of} $(f(M^{n-1}), x)$ if the dif\/ferential
$(f_t)_*$ has nullity $m$ at $x$.  In the case $\tilde{M}^n(c) = {\bf R}^n$, respectively~$S^n$,
the point $p=f_t(x)$ is a focal point
of multiplicity $m$ of $(f(M^{n-1}), x)$ if and only if $1/t$, respectively $\cot t$, is a principal curvature of multiplicity $m$ of $f(M^{n-1})$ at $x$ (see, for example, \cite[pp.~32--38]{Mil} or \cite[pp.~243--247]{CR}).

Let $f:M^{n-1} \rightarrow \tilde{M}^n(c)$ be an oriented hypersurface with constant
principal curvatures.  One can show from the formulas for the principal curvatures of a parallel hypersurface
(see, for example, \cite[pp.~131--132]{CR}) that if $f$ has constant principal curvatures, then each $f_t$ that is an
immersed hypersurface also has constant principal curvatures.  However, since the principal curvatures of $f$
are constant on $M^{n-1}$, the focal points along the normal geodesic to $f(M^{n-1})$ at $f(x)$ occur for the
same values of $t$ independent of the choice of point $x \in M^{n-1}$.  For example, for $\tilde{M}^n(c) = {\bf R}^n$,
if $\mu$ is a non-zero constant principal curvature of multiplicity $m>0$ of $M^{n-1}$, the map $f_{1/\mu}$ has constant rank
$n-1-m$ on $M^{n-1}$, and the set $f_{1/\mu}(M^{n-1})$ is an $(n-1-m)$-dimensional submanifold of ${\bf R}^n$ called
a {\em focal submanifold} of $f(M^{n-1})$.

One can show (see, for example, \cite[pp.~268--274]{CR}) that the level hypersurfaces of an isoparametric function $V$
form a family of parallel hypersurfaces (modulo reparametrization of the normal geodesics
to take into account the possibility that
$|{\rm grad}\, V|$ is not identically equal to one), and
each of these level hypersurfaces has constant principal curvatures.  Conversely,
one can begin with a connected hypersurface $f:M^{n-1} \rightarrow \tilde{M}^n(c)$ having constant principal
curvatures and construct an isoparametric function $V$ such that each
parallel hypersurface $f_t$ of $f$ is contained in a level set of $V$.
Therefore, one can def\/ine an {\em isoparametric hypersurface} to be a~hypersurface with constant principal curvatures,
and an isoparametric family of hypersurfaces can be characterized as a family of parallel hypersurfaces, each of which
has constant principal curvatures.  It is important that isoparametric hypersurfaces
always come as a family of parallel hypersurfaces together with their focal submanifolds.

As noted in the introduction, an isoparametric hypersurface in ${\bf R}^n$ must be an open
subset of a hyperplane, hypersphere or a spherical cylinder
$S^k \times {\mathbf R}^{n-k-1}$ (Somigliana \cite{Som} for $n=3$,
see also B.~Segre \cite{Seg1} and Levi-Civita \cite{Lev},
and B.~Segre \cite{Seg} for arbitrary $n$).  Shortly after the publication of the papers of Levi-Civita and Segre,
Cartan \cite{Car1,Car2,Car3,Car4} undertook the study of isoparametric hypersurfaces in arbitrary
real space-forms $\tilde{M}^n(c)$, $c \in {\bf R}$, and we now describe his primary contributions.

Let $f:M^{n-1} \rightarrow \tilde{M}^n(c)$ be an isoparametric hypersurface with $g$ distinct principal
curvatures $\mu_1,\ldots,\mu_g$, having respective multiplicities $m_1,\ldots,m_g$.
If $g>1$, Cartan showed that for each $i$,
$1 \leq i \leq g$,
\begin{gather}
\label{car-id}
\sum_{j\neq i} m_j \frac{c + \mu_i \mu_j}{\mu_i - \mu_j} = 0.
\end{gather}
This important equation, known as {\em Cartan's identity}, is crucial in Cartan's work on isoparametric hypersurfaces.
For example, using this identity, Cartan was able to classify isoparametric hypersurfaces in the cases $c \leq 0$
as follows.  In the case $c=0$, if $g=1$, then $f$ is totally umbilic, and it is well known that $f(M^{n-1})$
must be an open subset of a hyperplane or hypersphere.
If $g \geq 2$, then by taking an appropriate choice of unit normal f\/ield $\xi$, one
can assume that at least one of the principal curvatures is positive.  If $\mu_i$ is the smallest positive
principal curvature, then each term $\mu_i \mu_j/(\mu_i - \mu_j)$
in the sum in equation (\ref{car-id}) is non-positive, and thus must equal zero.
Therefore, there are at most two distinct principal curvatures, and if there are two, then one of them must be zero.
Hence, $g=2$ and one can show $f(M^{n-1})$ is an open subset of a spherical cylinder
by standard methods in Euclidean hypersurface theory.

\looseness=1
In the case $c=-1$, if $g=1$, then $f$ is totally umbilic, and it is well known that $f(M^{n-1})$
must be an open subset of a totally geodesic hyperplane, an equidistant hypersurface, a horosphere or a hypersphere
in $H^n$ (see, for example, \cite[p.~114]{Sp}).
If $g \geq 2$, then again one can arrange that at least one of the principal curvatures is positive.  Then
there must exist a positive principal curvature $\mu_i$ such that no principal curvature
lies between $\mu_i$ and $1/ \mu_i$. (Here $\mu_i$ is either the largest principal curvature between 0 and 1 or
the smallest principal curvature greater than or equal to one.) For this $\mu_i$,
each term $(-1 +\mu_i \mu_j)/(\mu_i - \mu_j)$
in the sum in equation (\ref{car-id}) is negative unless $\mu_j = 1/\mu_i$.  Thus, there are at most two distinct
principal curvatures, and if there are two, then they are reciprocals of each other.  Hence, $g=2$ and one can show
that $f(M^{n-1})$ is an open subset of a standard product
$S^k \times H^{n-k-1}$ in hyperbolic space $H^n$ (see Ryan~\cite{Ryan3}).

\looseness=1
In the sphere $S^n$, however, Cartan's identity does not lead to such strong restrictions on the number $g$
of distinct principal curvatures, and Cartan himself produced examples with $g = 1, 2, 3$ or 4 distinct
principal curvatures.  Moreover, he classif\/ied isoparametric hypersurfaces $f:M^{n-1} \rightarrow S^n$ with
$g \leq 3$ as follows.  In the case $g=1$, the hypersurface $f$ is totally umbilic, and it is well known that $f(M^{n-1})$
is an open subset of a great or small hypersphere in $S^n$ (see, for example, \cite[p.~112]{Sp}).
If $g=2$, then $M$ must be a standard product of two spheres,
\begin{gather}
\label{product}
S^p (r) \times S^q (s) \subset S^n (1) \subset {\bf R}^{p+1}
\times {\bf R}^{q+1} = {\bf R}^{n+1}, \qquad r^2 + s^2 = 1,
\end{gather}
where $n = p+q+1$ (see, for example, \cite[pp.~295--296]{CR}).

In the case of three distinct principal curvatures, Cartan \cite{Car3}
showed that all the principal curvatures must have the same multiplicity
$m=1,2,4$ or 8, and $f(M^{n-1})$ must be a tube of
constant radius over a standard embedding of a projective
plane ${\bf FP}^2$ into $S^{3m+1}$ (see, for example, \cite[pp.~296--299]{CR}),
where ${\bf F}$ is the division algebra
${\bf R}$, ${\bf C}$, ${\bf H}$ (quaternions),
${\bf O}$ (Cayley numbers), for $m=1,2,4,8,$ respectively. (In the case ${\bf F} = {\bf R}$, a standard embedding is a~Veronese surface in $S^4$.)
Thus, up to congruence, there is only one such family of isoparametric hypersurfaces
for each value of $m$.  For each of these hypersurfaces,
the focal set of $f(M^{n-1})$ consists of two antipodal standard embeddings of ${\bf FP}^2$,
and $f(M^{n-1})$ is a tube of constant radius over each focal submanifold.

In the process of proving this theorem, Cartan showed that any isoparametric family with $g$
distinct principal curvatures of the same multiplicity can be def\/ined by an equation of the form
$F = \cos gt$ (restricted to $S^n$), where $F$ is a harmonic homogeneous polynomial of degree $g$ on
${\bf R}^{n+1}$ satisfying
\begin{gather}
\label{car-eq}
|{\rm grad}\,F|^2 = g^2 r^{2g-2},
\end{gather}
where $r = |x|$ for $x \in {\bf R}^{n+1}$, and ${\rm grad}\, F$ is the gradient of $F$ in ${\bf R}^{n+1}$.  This was
a forerunner of M\"{u}nzner's general result that every isoparametric hypersurface is algebraic, and its
def\/ining polynomial satisf\/ies certain dif\/ferential equations which generalize those that Cartan found in this special case.

\looseness=1
In the case $g=4$, Cartan produced isoparametric hypersurfaces with four principal curvatures of
multiplicity one in $S^5$ and four principal curvatures of multiplicity two in $S^9$.  He noted all of his examples
are homogeneous, each being an orbit of a point under an appropriate closed subgroup of $SO(n+1)$.  Based on his results and the properties of his examples, Cartan asked the following three questions \cite{Car3}, all of which were answered
in the 1970's.

\begin{enumerate}\itemsep=0pt
\item  For each positive integer $g$, does there exist an isoparametric family with $g$ distinct principal
curvatures of the same multiplicity?

\item \looseness=-1 Does there exist an isoparametric family of hypersurfaces with more than three distinct~prin\-cipal curvatures such
that the principal curvatures do not all have the same multiplici\-ty?

\item\looseness=1 Does every isoparametric family of hypersurfaces admit a transitive group of isomet\-ries?
\end{enumerate}

Despite the depth and beauty of Cartan's work, the subject of isoparametric hypersurfaces in $S^n$
was virtually ignored for thirty years until a revival of the subject in the early 1970's by several authors.  Nomizu
\cite{Nom1,Nom2} wrote two papers describing the highlights of Cartan's work.  He also generalized
Cartan's example with four principal curvatures of multiplicity one to produce examples with four principal
curvatures having multiplicities $m_1 = m_2 = m$, and $m_3 = m_4 =1$, for any positive integer $m$.  This answered
Cartan's Question 2 in the af\/f\/irmative.  Nomizu also proved that every focal submanifold of every isoparametric
hypersurface must be a minimal submanifold of $S^n$.

Takagi and Takahashi \cite{TT} gave a complete classif\/ication of all homogeneous isoparametric hypersurfaces
in $S^n$, based on the work of Hsiang and Lawson \cite{HL}.  Takagi and Takahashi showed that each homogeneous
isoparametric hypersurface in $S^n$ is an orbit of the isotropy representation of a Riemannian symmetric space of
rank 2, and they gave a complete list of examples \cite[p.~480]{TT}.  This list contains examples with six
principal curvatures as well as those with $g=1,2,3,4$ principal curvatures, and in some cases with $g=4$, the principal
curvatures do not all have the same multiplicity, so this also provided an af\/f\/irmative answer to
Cartan's Question~2.

At about the same time as the papers of Nomizu and Takagi--Takahashi, M\"{u}nzner published two preprints that
greatly extended Cartan's work and have served as the basis for much of the research in the f\/ield since that
time.  The preprints were eventually published as papers \cite{Mu,Mu2} in 1980--1981.

In the f\/irst paper \cite{Mu} (see also Chapter~3 of \cite{CR}),
M\"{u}nzner began with a geometric study of the focal submanifolds of an isoparametric
hypersurface $f:M^{n-1} \rightarrow S^n$ with $g$ distinct principal curvatures.  Using the fact that each focal
submanifold of $f$ is obtained as a parallel map $f_t$, where $\cot t$ is a principal curvature of $f(M^{n-1})$,
M\"{u}nzner computed a formula for the shape operator of the focal submanifold $f_t(M^{n-1})$ in terms of the
shape operator of $f(M^{n-1})$ itself. In particular, his calculation shows that the assumption that $f(M^{n-1})$
has constant principal curvatures implies that the eigenvalues of the
shape operator $A_{\eta}$ are the same for every unit normal $\eta$ at every point of the focal submanifold
$f_t(M^{n-1})$.  Then by using a symmetry argument, he proved that if the principal curvatures
of $f(M^{n-1})$ are written as $\cot \theta_k, 0 < \theta_1 < \cdots < \theta_g < \pi$, with multiplicities $m_k$, then
\begin{gather}
\label{prin-curv}
\theta_k = \theta_1 + \frac{(k-1)}{g} \pi, \qquad 1 \leq k \leq g,
\end{gather}
and the multiplicities satisfy $m_k = m_{k+2}$ (subscripts mod $g$).  Thus, if $g$ is odd, all
of the multiplicities must be equal, and if $g$ is even, there are at most two distinct multiplicities.
M\"{u}nzner's calculation shows further that the focal submanifolds must
be minimal submanifolds of $S^n$, as Nomizu \cite{Nom1} had shown by a dif\/ferent proof, and that Cartan's
identity is equivalent to the minimality of the focal submanifolds (see also \cite[p.~251]{CR}).

If $\cot t$ is not a principal curvature of $f$, then the map $f_t$ in equation (\ref{parallel-S}) is also an
isoparametric hypersurface with $g$ distinct principal curvatures $\cot (\theta_1 - t),\ldots,\cot (\theta_g -t)$.
If $t = \theta_k$ (mod~$\pi$),
then the map $f_t$ is constant along each
leaf of the $m_k$-dimensional principal foliation~$T_k$, and the image of $f_t$ is a smooth focal submanifold of $f$ of
codimension $m_k +1$ in $S^n$. All of the hypersurfaces $f_t$ in a family of parallel isoparametric hypersurfaces
have the same focal submanifolds.

In a crucial step in the theory,
M\"{u}nzner then showed that $f(M^{n-1})$ and its parallel hypersurfaces and focal submanifolds are each contained in
a level set of a homogeneous polynomial~$F$ of degree $g$ satisfying the following
{\em Cartan--M\"{u}nzner differential equations}
on the Euclidean dif\/ferential operators ${\rm grad}\, F$ and Laplacian $\triangle F$ on ${\bf R}^{n+1}$,
\begin{gather}
\label{eq:C-M}
|{\rm grad}\,F|^2  =   g^2 r^{2g-2}, \qquad r = |x|, \\
 \triangle F  =  c r^{g-2}, \qquad c = g^2 (m_2-m_1)/2, \nonumber
\end{gather}
where $m_1$, $m_2$ are the two (possibly equal) multiplicities of the principal curvatures on $f(M^{n-1})$.
This generalized Cartan's polynomial in equation (\ref{car-eq}) for the case of~$g$ principal curvatures with
the same multiplicity.

Conversely, the level sets of the restriction $F|_{S^n}$ of a function $F$ satisfying equation (\ref{eq:C-M}) constitute
an isoparametric family of
hypersurfaces and their focal submanifolds, and $F$ is called the
{\em Cartan--M\"{u}nzner polynomial} associated to this family.
Furthermore, M\"{u}nzner showed that the level sets of $F$ are connected, and thus any
connected isoparametric hypersurface in $S^n$ lies in a unique compact, connected isoparametric hypersurface
obtained by taking the whole level set.

The values of the restriction $F|_{S^n}$
range between $-1$ and $+1$.  For $-1 < t < 1$, the level set $M_t = (F|_{S^n})^{-1}(t)$ is an isoparametric hypersurface,
while $M_{+} = (F|_{S^n})^{-1}(1)$ and $M_{-} = (F|_{S^n})^{-1}(-1)$ are focal submanifolds.  Thus, there are exactly
two focal submanifolds for the isoparametric family, regardless of the number $g$ of distinct principal curvatures.
Each principal curvature $\cot \theta_k$, $1 \leq k \leq g$, gives rise to two antipodal focal points corresponding
to the values $t = \theta_k$ and $t = \theta_k + \pi$ in equation (\ref{prin-curv}).  The $2g$ focal points are evenly
spaced at intervals of length $\pi/g$ along a normal geodesic to the isoparametric family, and they lie alternately
on the two focal submanifolds $M_{+}$ and $M_{-}$, which have respective codimensions $m_1 + 1$ and $m_2 + 1$ in~$S^n$.

Using this information, M\"{u}nzner showed that each isoparametric hypersurface $M_t$ in the family separates
the sphere $S^n$
into two connected components $D_1$ and $D_2$, such that $D_1$ is a~disk bundle with f\/ibers of dimension
$m_1 + 1$ over $M_{+}$,
and $D_2$ is a disk bundle with f\/ibers of dimension $m_2 + 1$ over $M_{-}$, where $m_1$ and $m_2$ are the multiplicities
of the principal curvatures that give rise to the focal submanifolds $M_{+}$ and $M_{-}$, respectively.

This topological situation has been the basis
for many important results in this f\/ield concerning the number $g$ of distinct principal curvatures and the multiplicities
$m_1$ and $m_2$. In particular, in his second paper, M\"{u}nzner \cite{Mu2} assumed that $M$ is a compact,
connected embedded hypersurface that separates $S^n$ into two disk bundles
$D_1$ and $D_2$ over compact manifolds with f\/ibers of dimensions $m_1 + 1$ and $m_2 + 1$, respectively.  From this
hypothesis, M\"{u}nzner proved that the dimension of the cohomological ring $H^*(M,{\bf Z}_2)$ must be $2\alpha$, where
$\alpha$ is one of the numbers  1, 2, 3, 4 or 6.  He then proved that if $M$ is a compact, connected isoparametric hypersurface
with $g$ distinct principal curvatures, then $\dim H^*(M,{\bf Z}_2) = 2g$.
Combining these two results, M\"{u}nzner obtained his major
theorem that the number $g$ of distinct principal curvatures of an isoparametric hypersurface in a sphere $S^n$
must be  1, 2, 3, 4 or 6.

M\"{u}nzner also obtained some restrictions on the possible values of the
multiplicities $m_1$ and~$m_2$.  These restrictions have been improved by several authors also using
topological arguments, as we will describe later in this section.  In particular, in the case $g=4$,
Abresch \cite{Ab} showed that if $m_1 = m_2 = m$,
then the only possible values for $m$ are 1 and 2, the values in
Cartan's examples.  In fact, up to congruence, Cartan's examples are the only isoparametric hypersurfaces with
four principal curvatures of the same multiplicity $m$.  This was proven by Takagi \cite{Takagi} for $m=1$ and
by Ozeki and Takeuchi \cite{OT} for $m=2$.

M\"{u}nzner's result gave a negative answer to Cartan's Question 1, but it pointed towards an af\/f\/irmative answer to
Cartan's Question 3, since the possible values  1, 2, 3, 4  or 6 for the number~$g$ of distinct principal curvatures
of an isoparametric hypersurface
agreed with the values of~$g$ for the homogeneous isoparametric
hypersurfaces on the list of Takagi and Takahashi.  However, in an important 2-part paper, Ozeki and
Takeuchi \cite{OT,OT2} used representations of Clif\/ford algebras to produce two inf\/inite series of
isoparametric families with four principal curvatures,
most of which are necessarily inhomogeneous, because their multiplicities do not
agree with the examples on the list of Takagi and Takahashi.  These papers of Ozeki and Takeuchi were then vastly
generalized by Ferus, Karcher and M\"{u}nzner \cite{FKM}, who also used representations of Clif\/ford algebras to
produce a large and important class of isoparametric families with four principal curvatures which contains
all known examples with $g=4$ with the exception of two homogeneous families.  We will now brief\/ly describe
the construction of Ferus, Karcher and M\"{u}nzner (see also, \cite[pp.~95--112]{Cec1}).

For each integer $m \geq 0$, the
{\em Clifford algebra} $C_m$ is the associative algebra over ${\bf R}$ that
is generated by a unity 1 and the elements $e_1,\ldots,e_m$ subject only to the relations
\begin{gather}
\label{cliff}
e_i^2 = -1, \qquad e_i e_j = - e_j e_i, \qquad i \neq j, \qquad 1 \leq i,j \leq m.
\end{gather}
One can show that the set
\begin{gather*}
\{1, e_{i_1} \cdots e_{i_r} \mid i_1 < \cdots < i_r, \  1 \leq r \leq m \},
\end{gather*}
forms a basis for the underlying vector space $C_m$, and thus $\dim C_m = 2^m$.

A representation of $C_m$ on ${\bf R}^l$ (of {\em degree} $l$) corresponds to a set of skew-symmetric matrices
$E_1,\ldots,E_m$ in the orthogonal group $O(l)$ such that
\begin{gather}
\label{clifford}
E_i^2 = -I, \qquad E_i E_j = - E_j E_i, \qquad i \neq j, \qquad 1 \leq i,j \leq m.
\end{gather}

Atiyah, Bott and Shapiro \cite{ABS} determined all of the Clif\/ford algebras according to the following table,
and they showed that the Clif\/ford algebra $C_{m-1}$ has an
irreducible representation of degree~$l$ if and only if $l = \delta (m)$ as in the table.
$$
\begin{array}{ccc}
{\underline m} \qquad & {\underline {C_{m-1}}} \qquad & {\underline {\delta(m)}} \\
 & & \\
1 \qquad & {\bf R} \qquad & 1\\
2 \qquad & {\bf C} \qquad & 2\\
3 \qquad & {\bf H} \qquad & 4\\
4 \qquad & {\bf H} \oplus {\bf H} \qquad & 4\\
5 \qquad & {\bf H}(2) \qquad & 8\\
6 \qquad & {\bf C}(4) \qquad & 8\\
7 \qquad & {\bf R}(8) \qquad & 8\\
8 \qquad & {\bf R}(8) \oplus {\bf R}(8) \qquad & 8\\
k+8 \qquad & C_{k-1} (16) \qquad & 16 \delta (k)
\end{array}
$$

\noindent
Reducible representations of $C_{m-1}$ on ${\bf R}^l$ for $l = k \delta (m)$, $k >1$, can be obtained
by taking a~direct sum of $k$ irreducible
representations of $C_{m-1}$ on ${\bf R}^{\delta(m)}$.

Given a representation of $C_{m-1}$ on ${\bf R}^l$ corresponding to skew-symmetric matrices $E_1,\ldots$, $E_{m-1}$
satisfying equation (\ref{clifford}),
Ferus, Karcher and M\"{u}nzner \cite{FKM} construct an isoparametric
family of hypersurfaces in the sphere $S^{2l-1} \subset {\bf R}^{2l} = {\bf R}^l \times {\bf R}^l$ with four
principal curvatures
such that one of the focal submanifolds is the manifold
\begin{gather}
\label{cliff-stief}
V_2 (C_{m-1}) = \Big\{(u,v) \in S^{2l-1} \mid |u|= |v| = \frac{1}{\sqrt{2}}, \ u \cdot v = 0, \ E_i u \cdot v = 0,
\ 1 \leq i \leq m-1 \Big\}.\!\!\!
\end{gather}
This is the {\em Clifford--Stiefel manifold}
of Clif\/ford orthogonal 2-frames of length
$1/\sqrt{2}$ in ${\bf R}^l$ (see Pinkall and Thorbergsson \cite{PT1}),
where vectors $u$ and $v$ in ${\bf R}^l$ are said to be
{\em Clifford orthogonal} if
$u \cdot v = E_1 u \cdot v = \cdots = E_{m-1} u \cdot v = 0$,
where $u \cdot v$ is the usual Euclidean inner product in ${\bf R}^l$.

Since $V_2 (C_{m-1})$ has codimension $m+1$ in $S^{2l-1}$, one of the principal curvatures of an isoparametric
hypersurface $M$ in this family has multiplicity $m_1 = m$.  The other multiplicity $m_2$ must satisfy the
equation $2m_1 + 2m_2 = \dim M = 2l-2$.  Therefore $m_2 = l-m-1$, where $l = k \delta (m)$ by the theorem of
Atiyah, Bott and Shapiro.

The isoparametric hypersurfaces resulting from this are called isoparametric hypersurfaces of {\em FKM-type}.
The following is a table of the multiplicities of the principal curvatures of the
FKM-hypersurfaces for multiplicities $(m_1,m_2) = (m, k \delta(m) - m - 1)$ for small values of $m$.
Of course, the multiplicity $m_2$ must be positive in order for this construction to lead to an isoparametric
hypersurface with four principal curvatures.  In the table below, the cases where $m_2 \leq 0$ are denoted by a dash.
$$
\begin{array}{llllllllllll}
\delta (m)| & 1 & 2 & 4 & 4 & 8 & 8 & 8 & 8 & 16 & 32 & \cdots\\
{\underline k} & & & & & & & & & & &\\
1 &-&-&-&-& (5,2) & (6,1) & - & - &(9,6) & (10,21)&  \cdots \\
2 &-&(2,1)&(3,4)&(4,3)&(5,10)&(6,9)&(7,8)&(8,7)&(9,22)&(10,53)&\cdots \\
3 &(1,1)&(2,3)&(3,8)&(4,7)&(5,18)&(6,17)&(7,16)&(8,15)&(9,38)&(10,85) & \cdots \\
4 &(1,2)&(2,5)&(3,12)&(4,11)&(5,26)&(6,25)&(7,24)&(8,23)&(9,54)& \cdot & \cdots\\
5 &(1,3)&(2,7)&(3,16)&(4,15)&(5,34)&(6,33)&(7,32)&(8,31)&\cdot & \cdot & \cdots \\
\cdot &\cdot&\cdot&\cdot&\cdot&\cdot&\cdot&\cdot&\cdot&\cdot&\cdot& \cdots \\
\cdot &\cdot&\cdot&\cdot&\cdot&\cdot&\cdot&\cdot&\cdot&\cdot&\cdot& \cdots \\
\cdot &\cdot&\cdot&\cdot&\cdot&\cdot&\cdot&\cdot&\cdot&\cdot&\cdot& \cdots \\
\end{array}
$$
\begin{center}
Multiplicities of the principal curvatures of FKM-hypersurfaces
\end{center}

\medskip

If $m \equiv 0$ (mod 4) and $l = k \delta (m)$, this construction yields $[k/2] + 1$ incongruent
isoparametric families having the same multiplicities.  Furthermore, the
families with multiplicities $(2,1)$, $(6,1)$, $(5,2)$ and one of the $(4,3)$-families are congruent to those
with multiplicities $(1,2)$, $(1,6)$, $(2,5)$ and $(3,4)$, respectively, and these are the only coincidences under
congruence among the FKM-hypersurfaces \cite{FKM}.

Many of the FKM examples are necessarily inhomogeneous, because their multiplicities do not agree with the multiplicities
of any hypersurface on the list of Takagi and Takahashi.  However,
Ferus, Karcher and M\"{u}nzner also gave a geometric proof of
the inhomogeneity of many of their examples through an examination
of the second fundamental forms of the focal submanifolds.
Later Q.-M. Wang \cite{Wang2} proved many results concerning the topology of the FKM examples
and Wu \cite{Wu2}
showed that there are only f\/initely many dif\/feomorphism classes
of compact isoparametric hypersurfaces with four
distinct principal curvatures.

All known examples of isoparametric
hypersurfaces with four principal curvatures are of FKM-type with the exception of two
homogeneous families, having multiplicities $(2,2)$ and~$(4,5)$. Beginning with M\"{u}nzner \cite{Mu,Mu2}, many mathematicians, including
Abresch \cite{Ab}, Grove and Halperin \cite{GH}, Tang \cite{Tang} and Fang \cite{Fang3,Fang1},
found restrictions on the multiplicities of the principal curvatures of an isoparametric hypersurface
with four principal curvatures. This series of results culminated with the paper of
Stolz \cite{Stolz}, who proved that the multiplicities of an isoparametric hypersurface with four principal curvatures
must be the same as those in the known examples of FKM-type or the two homogeneous exceptions.  All of these
papers are primarily topological in nature, based on M\"{u}nzner's result that an isoparametric hypersurface separates~$S^n$ into two disk bundles over the two focal submanifolds.  In particular, the proof of Stolz is homotopy theoretic,
and the main tools used are the Hopf invariant and the EHP-sequence.
In fact, Stolz proved his result under the more general assumption that $M$ is
a~compact, connected proper Dupin (not necessarily isoparametric) hypersurface embedded in~$S^n$.  Note that
Thorbergsson~\cite{Th1} had earlier shown that
a compact, connected proper Dupin hypersurface $M \subset S^n$ also separates
$S^n$ into two disk bundles over the f\/irst focal submanifolds on either side of $M$, as in
the isoparametric case.  This will be discussed
in more detail in the next section.

Cecil, Chi and Jensen \cite{CCJ1} then showed that if the multiplicities $(m_1,m_2)$ of
an isoparametric hypersurface $M \subset S^n$ with four principal curvatures satisfy $m_2 \geq 2 m_1 - 1$, then~$M$
must be of FKM-type.
Previously, Takagi \cite{Takagi} had shown that if one of the
multiplicities $m_1 =1$, then~$M$ is homogeneous and of FKM-type.
Ozeki and Takeuchi \cite{OT2} next showed that if $m_1 = 2$, then~$M$ is homogeneous
and of FKM-type,
except in the case of multiplicities $(2,2)$, in which case $M$ must be the homogeneous example of Cartan.
Taken together with these results of Takagi and Ozeki--Takeuchi, the theorem of Cecil, Chi and Jensen
classif\/ies isoparametric hypersurfaces with four principal curvatures for all
possible pairs of multiplicities except for four cases, the homogeneous pair  (4,5), and the FKM pairs
 (3,4), (6,9)  and  (7,8), for which the classif\/ication of isoparametric hypersurfaces remains an important open problem.

We now brief\/ly outline the proof of this result by Cecil, Chi and Jensen.
In Sections 8--9 of~\cite{CCJ1}, Cecil, Chi and Jensen use Cartan's method of moving frames to f\/ind necessary and suf\/f\/icient
conditions (equations (8.1)--(8.4) of~\cite{CCJ1}) for the codimension $m_1 + 1$
focal submanifold $M_+$ of an isoparametric hypersurface $M$ with four principal curvatures
and multiplicities $(m_1,m_2)$ to be a Clif\/ford--Stiefel manifold $V_2 (C_{m_1-1})$.
(Later Chi \cite{Chi} gave a dif\/ferent proof of the fact that equations (8.1)--(8.4) of~\cite{CCJ1} are necessary
and suf\/f\/icient to show that $M_+$ is a Clif\/ford--Stiefel manifold.)

These necessary and suf\/f\/icient
conditions involve the shape operators $A_{\eta}$ of $M_+$, where $\eta$ is a~unit normal vector to $M_+$ at a
point $x \in M_+$.  These shape operators are isospectral in that every~$A_{\eta}$ at every point $x \in M_+$
has the same eigenvalues $-1,0,1$, with respective multi\-pli\-ci\-ties~$m_2$, $m_1$, $m_2$.
If $\eta$ is a unit normal vector to $M_+$
at a point $x \in M_+$, then the point~$\eta$ is also in $M_+$ by M\"{u}nzner's results,
since it lies at a distance $\pi/2$ along the
normal geodesic to $M_+$ beginning at the point $x$ in the direction $\eta$.  The shape operators corresponding
to an orthonormal basis of normal vectors to $M_+$ at the point $x$
determine a family of $m_1 + 1$ homogeneous polynomials. Similarly, the shape operators corresponding
to an orthonormal basis of normal vectors to $M_+$ at the point $\eta$ determine a family of $m_1 + 1$
homogeneous polyno\-mials.
In Section 10 of \cite{CCJ1}, Cecil, Chi and Jensen show that
these two families of polynomials have the same zero set in projective space
by use of a formulation of the Cartan--M\"{u}nzner polynomial due to Ozeki and Takeuchi \cite{OT}.
Finally, in Sections 11--13 of \cite{CCJ1}, the authors employ
techniques from algebraic geometry to show that the fact that these two sets of polynomials have the same zero set
leads to a proof that the necessary and
suf\/f\/icient conditions for $M_+$ to be a Clif\/ford--Stiefel manifold are satisf\/ied if $m_2 \geq 2 m_1 - 1$.
This completes the proof that $M$ is of FKM-type, since $M$ is a tube of constant radius over
the Clif\/ford--Stiefel mani\-fold~$M_+$.

After this, Immervoll \cite{Im} gave a dif\/ferent proof of the theorem of Cecil, Chi and Jensen using isoparametric
triple systems.  The use of triple systems to study isoparametric hypersurfaces was introduced in
a series of papers in the 1980's by Dorfmeister and Neher \cite{DN1,DN2,DN3,DN4,DN,DN6}.

In the case of an isoparametric hypersurface with six principal curvatures, M\"{u}nzner showed
that all of the principal curvatures must have the same multiplicity $m$, and
Abresch \cite{Ab} showed that $m$ must equal 1 or 2.
By the classif\/ication of homogeneous isoparametric hypersurfaces due to Takagi and Takahashi \cite{TT},
there is only one homogeneous family in each case up to congruence.  In the case of multiplicity $m=1$,
Dorfmeister and Neher \cite{DN} showed that an isoparametric
hypersurface must be homogeneous, thereby completely classifying that case.
The proof of Dorfmeister and Neher is quite algebraic in nature, and recently Miyaoka \cite{Mi10} has
given a shorter, more geometric proof of this result.

Miyaoka \cite{Mi6} also
gave a geometric description of the case $m=1$, showing that a
homogeneous isoparametric hypersurface $M^6$ in $S^7$
can be obtained as the inverse image under
the Hopf f\/ibration $h:S^7 \rightarrow S^4$ of
an isoparametric hypersurface with three principal curvatures of
multiplicity one in $S^4$. Miyaoka also showed that the two focal
submanifolds of $M^6$ are not congruent, even though they
are lifts under $h^{-1}$ of congruent Veronese surfaces in $S^4$.
Thus, these focal submanifolds are two
non-congruent minimal homogeneous embeddings of
${\bf RP}^2 \times S^3$ in $S^7$.

Peng and Hou \cite{Peng}
gave explicit forms for the Cartan--M\"{u}nzner polynomials of degree six
for the homogeneous isoparametric hypersurfaces with $g=6$, and
Fang \cite{Fang2}
proved several results concerning the topology of isoparametric
and compact proper Dupin hypersurfaces with six principal
curvatures.

The classif\/ication of isoparametric hypersurfaces with six principal curvatures of multiplicity $m=2$ is part of
Problem 34 on Yau's \cite{Yau} list of important open problems in geometry, and it
remains an open problem.
It has long been conjectured that the one homogeneous family in the case $g=6$, $m=2$,
is the only isoparametric family in this case,
but this conjecture has resisted proof for a long time.  The approach that Miyaoka
\cite{Mi10} used in the case $m=1$ shows promise of possibly leading to a proof of this conjecture,
but so far a complete proof has not been published.

In the 1980's, the notion of an isoparametric hypersurface was extended to submanifolds of
codimension greater than one in $S^n$ by several authors independently. (See
Carter and West \cite{CW1,West}, Harle~\cite{Har},
Str\"{u}bing~\cite{Str} and Terng~\cite{Te1}.)  An immersed
submanifold $\phi:V \rightarrow {\bf R}^n$ (or~$S^n$) is def\/ined to be
{\em isoparametric} if its normal bundle $N(V)$ is f\/lat,
and if for any locally def\/ined normal f\/ield $\xi$ which is parallel
with respect to the normal connection~$\nabla^{\perp}$,
the eigenvalues of the shape operator $A_{\xi}$ are constant.
This theory was then developed extensively by many authors over the next decade,
especially in the papers of Terng \cite{Te1}, Palais and Terng~\cite{PalT1},
and Hsiang, Palais and Terng~\cite{HPT} (see also the book \cite{PalT}).
Finally, Thorbergsson~\cite{Th2} used the theory of Tits buildings to
show that all irreducible
isoparametric submanifolds of codimension greater than one in $S^n$ are homogeneous, and therefore they are
principal orbits of isotropy representations of symmetric spaces,
also known as generalized f\/lag manifolds or standard embeddings of $R$-spaces.
(See Bott and Samelson \cite{BS}, Takeuchi and Kobayashi~\cite{TK}, Dadok~\cite{Dad}, Hahn~\cite{Hahn2}
or the book by Berndt, Console and Olmos~\cite{BCO}.)
Later Olmos \cite{Olm} and Heintze and Liu \cite{HLiu} gave alternate proofs of Thorbergsson's result.

Heintze, Olmos and Thorbergsson \cite{HOT} def\/ined a submanifold $\phi:V \rightarrow {\bf R}^n$  (or $S^n$)
to have {\em constant principal curvatures} if for any smooth
curve $\gamma$ on $V$ and any parallel normal vector f\/ield~$\xi(t)$ along $\gamma$, the shape operator $A_{\xi(t)}$ has
constant eigenvalues along~$\gamma$.  If the normal bund\-le~$N(M)$ is f\/lat, then having constant principal curvatures is
equivalent to being isoparametric.  They then
showed that a submanifold with constant principal curvatures is
either isoparametric or a focal submanifold of an isoparametric
submanifold. The excellent survey article of Thorbergsson \cite{Th6} gives a more detailed account of the theory of
isoparametric submanifolds of codimension greater than one in~$S^n$.

An immersion $\phi:V \rightarrow {\bf R}^n$ of a compact, connected manifold $V$ into ${\bf R}^n$
is said to be {\em taut} (see, for example, \cite{CW} or \cite{CR}) if every Morse
function of the form $L_p:V \rightarrow {\bf R}$, where $L_p(x) = |p - \phi (x)|^2$,
for $p \in {\bf R}^n$, has the minimum
number of critical points required by the Morse inequalities using ${\bf Z}_2$-homology, i.e., it is a perfect
Morse function.
Tautness can also be studied for submanifolds
of $S^n$ using spherical distance functions instead of Euclidean distance functions, and tautness is invariant under stereographic projection
and its inverse map.

Isoparametric submanifolds (of any codimension) and their focal submanifolds are all taut submanifolds of $S^n$.
This was proven by Cecil and Ryan \cite{CRBer} for hypersurfaces and by Hsiang, Palais and Terng \cite{HPT} for
isoparametric submanifolds of codimension greater
than one in $S^n$, and this was an important fact in the general development of the theory.

Carter and West \cite{CW3,CW4,CW6,CW5},
introduced the notion of totally focal submanifolds and studied its relationship to the
isoparametric property. A submanifold $\phi:V \rightarrow {\bf R}^n$ is said to be
{\em totally focal} if the critical points of
every Euclidean distance function $L_p$ on $V$
are either all non-degenerate or all degenerate.  An isoparametric submanifold in $S^n$ is
totally focal, and the main result of \cite{CW5} is that a
totally focal submanifold must be isoparametric.  However,
Terng and Thorbergsson \cite[p.~197]{TTh1} have noted that there is
a gap in the proof of this assertion, specif\/ically in the proof
of Theorem~5.1 of \cite{CW5}.

\looseness=1
Wu \cite{Wu} and Zhao \cite{Zhao} generalized the theory of isoparametric submanifolds of codimension greater than one
to submanifolds of hyperbolic space, and  Verh\'{o}czki \cite{Ver} developed a theory of
isoparametric submanifolds for Riemannian manifolds
which do not have constant curvature.
West \cite{West1} and Mullen \cite{Mul} formulated a theory of isoparametric
systems on symmetric spaces, and Terng and Thorbergsson \cite{TTh}
studied compact isoparametric submanifolds of symmetric spaces
using the related notion of equifocal submanifolds.  Tang \cite{Tang2} then did a thorough study of the possible
multiplicities of the focal points of equifocal hypersurfaces in symmetric spaces.  In a generalization
of Thorbergsson's result for submanifolds of $S^n$,
Christ~\cite{Christ} proved that a complete connected irreducible equifocal submanifold of codimension greater than one in a~simply connected compact symmetric space is homogeneous. Finally,
a promising recent generalization of the theory of isoparametric submanifolds is the theory of singular Riemannian foliations
admitting sections (see Alexandrino~\cite{Alex}, T\"{o}ben~\cite{Toben}, and Lytchak and Thorbergsson~\cite{LT}).

In a later paper, Terng and Thorbergsson \cite{TTh1}
gave a def\/inition of tautness for submanifolds of arbitrary complete
Riemannian manifolds, and they discussed the notions of isoparametric, equifocal and Dupin submanifolds in that setting.
In a related development, Carter and \c{S}ent\"{u}rk~\cite{CS} studied the space of immersions parallel to a given immersion whose
normal bundle has trivial holonomy group.

Terng \cite{Te5} considered isoparametric submanifolds in inf\/inite-dimensional
Hilbert spaces and generalized many results from the f\/inite-dimensional case to that setting.  Pinkall
and Thorbergsson \cite{PT3} then gave more examples of such submanifolds, and
Heintze and Liu \cite{HLiu} generalized the f\/inite-dimensional homogeneity result of Thorbergsson~\cite{Th2}
to the inf\/inite-dimensional case.

\looseness=1
Nomizu \cite{Nom3} began the study of isoparametric hypersurfaces in pseudo-Riemannian space forms
by proving a generalization of Cartan's identity for space-like hypersurfaces in a Lorentzian space form
$\tilde{M}^n_1(c)$ of constant sectional curvature $c$.  As a consequence of this identity, Nomizu showed that
a space-like isoparametric hypersurface in $\tilde{M}^n_1(c)$ can have at most two distinct principal curvatures
if $c \geq 0$. Recently, Li and Xie \cite{LX} have shown that this conclusion also holds for space-like isoparametric
hypersurfaces in $\tilde{M}^n_1(c)$ for $c<0$.
Magid \cite{Mag} studied isoparametric hypersurfaces in Lorentz space whose shape operator is not diagonalizable,
and Hahn~\cite{Hahn} contributed an extensive study of isoparametric hypersurfaces in pseudo-Riemannian space forms
of arbitrary signatures.
Recently Geatti and Gorodski~\cite{Geatti} have extended this theory
further by showing that a polar orthogonal representation of a connected real reductive
algebraic group has the same closed orbits as the isotropy representation of a pseudo-Riemannian symmetric space.
Working in a dif\/ferent direction, Niebergall and Ryan \cite{NieR1,NieR2,NieR3,NieR4} developed a~theory
of isoparametric and Dupin hypersurfaces in af\/f\/ine dif\/ferential geo\-metry.

There is also an extensive theory of real hypersurfaces with constant principal curvatures
in complex space forms that is closely related to the theory
of isoparametric hypersurfaces in spheres.  See the survey articles of Niebergall and Ryan \cite{NieR}
and Berndt \cite{Berndt} for more detail.

In applications to Riemannian geometry, Solomon \cite{Sol1,Sol2,Sol3} found results concerning the spectrum
of the Laplacian of isoparametric hypersurfaces in $S^n$ with
three or four principal curvatures, and
Eschenburg and Schroeder \cite{ES} studied the behavior of the Tits metric on isoparametric
hypersurfaces.  Finally, Ferapontov \cite{Fera,Fera1}
studied the relationship between isoparametric and
Dupin hypersurfaces and Hamiltonian
systems of hydrodynamic type and listed several open research
problems in that context.

\section{Dupin hypersurfaces}
\label{sec:2}
As noted in the introduction, the theory of proper Dupin hypersurfaces is closely related to that of isoparametric
hypersurfaces.  In contrast to the situation for isoparametric hypersurfaces, however, there are both local
and global aspects to the theory of proper Dupin hypersurfaces with quite dif\/ferent results concerning the number
of distinct principal curvatures and their multiplicities, for example.  In this section, we survey the primary results
in this f\/ield.  For the sake of completeness, we will begin with a formal def\/inition of the Dupin
and proper Dupin properties.

Via stereographic projection $\tau :S^n - \{P\} \rightarrow {\bf R}^n$ with pole $P \in S^n$
and its inverse map $\sigma$ from ${\bf R}^n$ into $S^n$, the theory of
Dupin hypersurfaces is essentially the same for hypersurfaces in ${\bf R}^n$ or $S^n$
(see, for example, \cite[pp.~132--151]{CR}),
and we will use whichever ambient space is most convenient for the discussion at hand.  Since we have been dealing
with isoparametric hypersurfaces in spheres, we will formulate our def\/initions here in terms of hypersurfaces in~$S^n$.

Let $f:M \rightarrow S^n$
be an immersed hypersurface, and let $\xi$ be a locally def\/ined
f\/ield of unit normals to $f(M)$.
A {\em curvature surface} of $M$ is a smooth submanifold $S$
such that for each point $x \in S$, the tangent space
$T_xS$ is equal to a principal space of the shape operator
$A$ of $M$ at $x$. The hypersurface $M$ is said to be {\em Dupin} if:

\begin{enumerate}\itemsep=0pt
\item[(a)] along each curvature surface, the corresponding principal
curvature is constant.
\end{enumerate}
The hypersurface $M$ is called {\em proper Dupin} if, in addition
to Condition~(a), the following condition is satisf\/ied:

\begin{enumerate}\itemsep=0pt
\item[(b)] the number $g$ of distinct principal curvatures is constant on
$M$.
\end{enumerate}
An obvious and important class of proper Dupin hypersurfaces are the isoparametric hypersurfaces in
$S^n$, and those hypersurfaces in ${\bf R}^n$ obtained from isoparametric hypersurfaces in $S^n$ via
stereographic projection.  For example, the well-known ring cyclides of Dupin in ${\bf R}^3$ are obtained
from a standard product torus $S^1(r) \times S^1(s) \subset S^3$, $r^2+s^2=1$, in this way.  These examples will
be discussed in more detail later in this section.

\looseness=1
We begin by mentioning several well-known basic facts about Dupin hypersurfaces which are proven in Section 2.4 of \cite{CR},
for example.  As above, let $f:M \rightarrow S^n$
be an immersed hypersurface, and let $\xi$ be a locally def\/ined
f\/ield of unit normals to $f(M)$.
First by the Codazzi equation, Condition (a) is automatically satisf\/ied on a curvature surface $S$
of dimension greater than one.
Second, Condition (b) is equivalent to requiring that each continuous
principal curvature have constant multiplicity on $M$. Further, the number of distinct principal
curvatures is locally constant on a dense open subset of any
hypersurface in $S^n$ (see Singley~\cite{Sin}).

Next, if a continuous principal curvature function $\mu$ has constant
multiplicity $m$ on a connected open subset $U \subset M$, then $\mu$ is a smooth function, and
the distribution $T_{\mu}$ of principal
vectors corresponding to $\mu$ is a smooth foliation whose leaves are the curvature surfaces corresponding to $\mu$. The principal curvature
$\mu$ is constant along each of its curvature surfaces in $U$
if and only if these curvature surfaces are open subsets
of $m$-dimensional great or small spheres in $S^n$.  Suppose that $\mu = \cot \theta$, for $0 < \theta < \pi$,
where $\theta$ is a smooth function on $U$.
The corresponding
{\em focal map} $f_{\mu}$
which maps $x\in M$ to the focal point $f_{\mu}(x)$ is given by the formula,
\begin{gather}
\label{focal map}
f_{\mu}(x) = \cos \theta (x) \; f(x) + \sin \theta (x) \; \xi(x).
\end{gather}
The principal curvature $\mu$ also determines a second focal map, whose image is antipodal to the image of $f_{\mu}$,
obtained by replacing $\theta$ by $\theta + \pi$ in
equation (\ref{focal map}).
The principal curvature function~$\mu$ is constant along each of its curvature surfaces in $U$ if and only if
the focal map~$f_{\mu}$ factors through an immersion of the
$(n-1-m)$-dimensional space of leaves $U/T_{\mu}$ into~$S^n$, and thus~$f_{\mu}(U)$ is an $(n-1-m)$-dimensional
submanifold of $S^n$.

\looseness=1
The {\em curvature sphere} $K_{\mu}(x)$ corresponding to
the principal curvature $\mu$ at a point $x \in U$ is the hypersphere in $S^n$ through $f(x)$ centered at the focal
point $f_{\mu}(x)$.  Thus, $K_{\mu}(x)$ is tangent to $f(M)$ at $f(x)$.  The principal curvature $\mu$ is constant
along each of its curvature surfaces if and only if the curvature sphere map $K_{\mu}$ is constant
along each of these curvature surfaces.

Thus, on an open subset $U$ on which Condition (b) holds, Condition (a) is equivalent
to requiring that each curvature surface in each principal
foliation be an open subset of a metric sphere in $S^n$
of dimension equal to the multiplicity of the corresponding
principal curvature. Condition (a) is also equivalent to the condition that
along each curvature surface, the corresponding curvature sphere map is constant. Finally, on $U$, Condition (a)
is equivalent to requiring that for each principal curvature $\mu$, the image of the focal
map $f_{\mu}$ is a smooth submanifold of $S^n$ of codimension $m+1$, where $m$ is the multiplicity of $\mu$.

\looseness=1
One consequence of these results
is that like isoparametric hypersurfaces, proper Duper hypersurfaces are algebraic.  For simplicity, we will formulate
this result in terms of hypersurfaces of ${\bf R}^n$.  The theorem states that a connected proper Dupin hypersurface
$f:M \rightarrow{\bf R}^n$ must be contained in a connected component of an
irreducible algebraic subset of ${\bf R}^n$ of dimension $n-1$.  This result was widely thought to be true in the 1980's,
and indeed Pinkall \cite{P6} sent the author a letter in 1984 that contained a sketch of a proof of this result.
However, a proof was not published until recently by Cecil, Chi and Jensen \cite{CCJ3}, who used methods of real
algebraic geometry to give a complete proof based on Pinkall's sketch. The proof makes use of the various principal foliations whose leaves are open subsets of spheres to construct an analytic algebraic parametrization of a
neighborhood of $f(x)$ for each point $x \in M$.
In contrast to the situation for isoparametric hypersurfaces, however, a connected proper Dupin hypersurface
in $S^n$ does not necessarily lie in a compact connected proper Dupin hypersurface, as discussed later in this section.
The algebraicity,
and hence analyticity, of proper Dupin hypersurfaces was useful in clarifying certain f\/ine points in the
recent paper \cite{CCJ3} of Cecil, Chi and Jensen on proper Dupin hypersurfaces with four principal curvatures.

The def\/inition of Dupin can be extended to submanifolds of
codimension greater than one as follows.  Let
$\phi:V \rightarrow {\bf R}^n$ (or $S^n$) be a submanifold of
codimension greater than one, and let $B^{n-1}$ denote the
unit normal bundle of $\phi(V)$.  A {\em curvature
surface} (see Reckziegel \cite{Reck}) is a connected submanifold $S \subset V$ for which
there is a parallel section $\eta : S \rightarrow B^{n-1}$
such that for each $x \in S$, the tangent space $T_xS$ is
equal to some smooth eigenspace of the shape operator
$A_{\eta}$. The submanifold $\phi(V)$ is said to be {\em Dupin}
if along each curvature surface, the corresponding principal
curvature of $A_{\eta}$ is constant, and a Dupin submanifold is {\em proper Dupin} if the number
of distinct principal curvatures is constant
on the unit normal bundle $B^{n-1}$.  An isoparametric submanifold of codimension greater than one is always
Dupin, but it may not be proper Dupin.  (See \cite[pp.~464--469]{Te2} for more detail on this point.)
Pinkall \cite[p.~439]{P2} proved that every extrinsically
symmetric submanifold of a real space form is Dupin.  Takeuchi~\cite{Tak} then determined which of these are proper Dupin.

\looseness=1
Pinkall \cite{P2,P3} situated the study of Dupin hypersurfaces in the context of Lie sphere geometry, and
this approach has proven to be very useful in subsequent research in the f\/ield.  In particular,
Dupin submanifolds of codimension greater than one in $S^n$ can be studied naturally in this setting.
Here we give a brief outline of this
approach to Dupin hypersurfaces in $S^n$.  (See Pinkall \cite{P3}, Chern \cite{Chern},
Cecil and Chern \cite{CC1}, or the book \cite{Cec1} for more detail.)

A {\em Lie sphere} is an oriented hypersphere or a point sphere (zero radius) in $S^n$.
The set of all Lie spheres in $S^n$
is in bijective correspondence with the set of points $[x] = [(x_1,\ldots,x_{n+3})]$ (homogeneous coordinates)
in real projective space ${\bf P}^{n+2}$ that lie
on the quadric hypersurface~$Q^{n+1}$ (the {\em Lie quadric})
determined by the equation
$\langle x, x \rangle = 0$, where
\begin{gather}
\label{Lie-metric}
\langle x, y \rangle = -x_1 y_1 + x_2 y_2 + \cdots +x_{n+2} y_{n+2} - x_{n+3} y_{n+3}
\end{gather}
is a bilinear form of signature $(n+1,2)$ on the indef\/inite inner product space ${\bf R}^{n+3}_2$. Here the sphere $S^n$ is
identif\/ied with the unit sphere in ${\bf R}^{n+1} \subset {\bf R}^{n+3}_2$, where ${\bf R}^{n+1}$ is spanned by the
standard basis vectors $\{e_2,\ldots,e_{n+2}\}$ in ${\bf R}^{n+3}_2$.
Specif\/ically, the oriented hypersphere with center $p \in S^n$ and
signed radius $\rho$ (the sign denotes the orientation) corresponds to the point
$[(\cos \rho , p, \sin \rho )]$ in $Q^{n+1}$. Point
spheres $p$ in $S^n$ correspond to those points $[(1, p, 0)]$ in $Q^{n+1}$ with radius $\rho = 0$.

The Lie quadric $Q^{n+1}$ contains projective lines but no linear
subspaces of ${\bf P}^{n+2}$ of higher dimension.  If $[x]$ and $[y]$ are points in~$Q^{n+1}$, the line
$[x,y]$ determined by $[x]$ and $[y]$
lies on~$Q^{n+1}$ if and only if $\langle x,y \rangle=0$.
This condition means that the hyperspheres in $S^n$ corresponding
to~$[x]$ and~$[y]$ are in oriented contact, i.e., they are tangent and have the same orientation at the point of contact.
For a point sphere and an oriented sphere, oriented contact means that the point lies on the sphere.

A {\em Lie sphere transformation} is a projective
transformation of ${\bf P}^{n+2}$ that maps $Q^{n+1}$ to
itself.  In terms of the geometry of $S^n$, a Lie sphere
transformation is a dif\/feomorphism on the space of Lie spheres that preserves oriented contact
of spheres (see Pinkall \cite[p.~431]{P3}), since it takes lines on $Q^{n+1}$ to lines on
$Q^{n+1}$.  The group of Lie sphere transformations is
isomorphic to $O(n+1,2)/ \{ \pm I \}$, where $O(n+1,2)$
is the orthogonal group for the metric in equation (\ref{Lie-metric}).
A~Lie sphere transformation that takes point spheres to point spheres is a {\em M\"{o}bius transformation}, i.e.,
it is induced by a conformal dif\/feomorphism of~$S^n$, and the set of all M\"{o}bius transformations is a subgroup of the
Lie sphere group.
An example of a Lie sphere transformation that is not
a~M\"{o}bius transformation is a parallel transformations~$P_t$, which f\/ixes the center of each Lie
sphere but adds $t$ to its signed radius.
The group of Lie sphere transformations is generated by the union of the M\"{o}bius group and the group of all
parallel transformations.

The manifold  $\Lambda^{2n-1}$ of projective lines on $Q^{n+1}$ has a {\em contact structure},
i.e., a globally def\/ined 1-form $\omega$ such that
$\omega \wedge d\omega ^{n-1}$ never vanishes on
$\Lambda ^{2n-1}$.  The condition $\omega = 0$ def\/ines a~codimension one distribution $D$ on $\Lambda ^{2n-1}$ which has
integral submanifolds of dimension $n-1$ but none of higher
dimension.  A {\em Legendre submanifold} is one of these integral
submanifolds of maximal dimension, i.e.,
an immersion $\lambda : M^{n-1} \rightarrow
\Lambda ^{2n-1}$ such that $\lambda ^* \omega = 0$.

An oriented hypersurface $f:M^{n-1} \rightarrow S^n$ with f\/ield of unit
normals $\xi :M^{n-1} \rightarrow S^n$ naturally induces
a Legendre submanifold $\lambda = [k_1, k_2]$ (the line determined by the points $[k_1]$ and $[k_2]$ in~$Q^{n+1}$), where
\begin{gather}
k_1 = (1,f,0), \qquad k_2 = (0,\xi ,1).
\end{gather}
For each $x \in M^{n-1}$, $[k_1(x)]$ is the unique point sphere and
$[k_2 (x)]$ is the unique great sphere in the parabolic pencil of spheres in $S^n$ corresponding to
the points on the line $\lambda (x)$.
Similarly, an immersed submanifold $\phi :V \rightarrow S^n$ of codimension
greater than one induces a Legendre submanifold
whose domain is the bundle $B^{n-1}$ of unit normal
vectors to $\phi (V)$.  In each case, $\lambda$ is called the {\em Legendre lift} of the submanifold in $S^n$.
Note that the space of Legendre submanifolds is larger than the space of Legendre lifts of
submanifolds of~$S^n$, since the point sphere map of an arbitrary Legendre submanifold may not have constant rank
as a map into~$S^n$.

A Lie sphere transformation $\beta$ maps lines on
$Q^{n+1}$ to lines on $Q^{n+1}$, so it naturally
induces a map $\tilde{\beta }$ from
$\Lambda ^{2n-1}$ to itself.  If $\lambda$ is a
Legendre submanifold, then $\tilde{\beta }\lambda$
is also a Legendre submanifold, which is denoted
$\beta \lambda$ for short.  These two Legendre
submanifolds are said to be {\em Lie equivalent}.
We will also say that two submanifolds of $S^n$ are Lie equivalent, if their corresponding Legendre lifts
are Lie equivalent.
If $\beta$ is a M\"{o}bius transformation, then the two Legendre
submanifolds are said to be {\em M\"{o}bius equivalent}.
Finally, if $\beta$ is the parallel transformation $P_t$ and~$\lambda$ is the Legendre submanifold induced by an
oriented hypersurface $f:M \rightarrow S^n$, then
$P_t\lambda$ is the Legendre lift of the
parallel hypersurface $f_{-t}$ at oriented distance $-t$ from $f$ (see, for example, \cite[p.~67]{Cec1}).

To def\/ine the Dupin property for a Legendre submanifold
$\lambda : M^{n-1} \rightarrow \Lambda ^{2n-1}$, one replaces
the principal curvature function on $M^{n-1}$ in Conditions (a) and (b) above
with the corresponding curvature sphere map $K:M^{n-1} \rightarrow Q^{n+1}$,
which is naturally def\/ined in this setting.
One can easily show that a Lie sphere transformation $\beta$ maps curvatures spheres of $\lambda$ to curvature
spheres of $\beta \lambda$, and that Conditions~(a) and~(b) are preserved by $\beta$ (see \cite[pp.~67--70]{Cec1}).
Thus, both the Dupin and proper Dupin properties are invariant under Lie sphere transformations.

Pinkall \cite{P3} introduced four standard constructions for obtaining a proper Dupin hypersurface in
${\bf R}^{n+m}$ from a proper Dupin hypersurface in ${\bf R}^n$.  We f\/irst describe
Pinkall's constructions in the case $m=1$.

Start with a Dupin hypersurface
$W^{n-1}$ in ${\bf R}^n$ and then consider ${\bf R}^n$
as the linear subspace ${\bf R}^n \times \{ 0 \}$
in ${\bf R}^{n+1}$.  The following
constructions yield a Dupin hypersurface $M^n$ in
${\bf R}^{n+1}$.
\begin{enumerate}\itemsep=0pt
\item[(1)] Let $M^n$ be the cylinder $W^{n-1} \times {\bf R}$ in
${\bf R}^{n+1}$.
\item[(2)] Let $M^n$ be the hypersurface in ${\bf R}^{n+1}$
obtained by rotating
$W^{n-1}$ around an axis
${\bf R}^{n-1} \subset {\bf R}^n$.
\item[(3)] Project $W^{n-1}$ stereographically onto a hypersurface
$V^{n-1} \subset S^n \subset {\bf R}^{n+1}$.  Let
$M^n$ be the cone over $V^{n-1}$ in ${\bf R}^{n+1}$.
\item[(4)] Let $M^n$ be a tube in ${\bf R}^{n+1}$ around $W^{n-1}$.
\end{enumerate}
Even though Pinkall gave these four constructions, he noted that the tube construction and the cone construction are
Lie equivalent \cite[p.~438]{P3}, and therefore in the context of Lie sphere geometry, it is suf\/f\/icient
to deal only with the tube, cylinder and surface of revolution constructions.

In general, these constructions introduce a new principal curvature
of multiplicity one which is easily seen to be constant along its lines
of curvature.  The other principal curvatures are determined by the
principal curvatures of $W^{n-1}$, and the Dupin property is preserved
for these principal curvatures (see \cite{CecGD} or \cite[pp.~127--141]{Cec1} for a full discussion of these constructions
in the setting of Lie sphere geometry). Thus,
if $W^{n-1}$ is a proper Dupin hypersurface in ${\bf R}^n$
with~$g$ distinct principal curvatures, then in general, $M^n$ is a
proper Dupin hypersurface in ${\bf R}^{n+1}$ with~$g+1$ distinct principal curvatures.
However, this is not always the case, as we will explain below.
Pinkall also pointed out that these constructions can easily be generalized to produce a new principal curvature
of multiplicity $m > 1$ by considering the constructions in ${\bf R}^n \times {\bf R}^m$ instead of
${\bf R}^n \times {\bf R}$.

Let us now examine the possible problems encountered in trying to obtain a proper Dupin hypersurface by
each of the three constructions. In all cases, suppose that $W^{n-1}$ is a proper Dupin hypersurface in ${\bf R}^n$
with $g$ distinct principal curvatures.

In the cylinder construction,
the new principal curvature of $M^n$ is identically zero while
the other principal curvatures of $M^n$ are equal to those
of $W^{n-1}$.  Thus if one of the principal curvatures $\kappa$
of $W^{n-1}$ is zero at some points but not identically
zero, then the number of distinct principal curvatures
is not constant on $M^n$, and $M^n$ is Dupin but not proper Dupin.

If $M^n$ is a tube in ${\bf R}^{n+1}$ of radius
$\epsilon$ over $W^{n-1}$, then there are exactly two
distinct principal curvatures at the points in the set
$W^{n-1} \times \{ \pm \epsilon \}$ in $M^n$, regardless
of the number of distinct principal curvatures on $W^{n-1}$.
Thus, $M^n$ is not a proper Dupin hypersurface unless the original hypersurface $W^{n-1}$ is totally umbilic, i.e., it
has only one distinct principal curvature at each point.

For the surface of revolution construction, if the focal point corresponding to a principal curvature
$\kappa$ at a point $x$ of the prof\/ile submanifold $W^{n-1}$ lies on the axis of revolution ${\bf R}^{n-1}$, then
the principal curvature of $M^n$ at $x$
determined by $\kappa$ is equal to the new principal curvature of $M^n$ resulting from
the surface of revolution construction.
Thus, if the focal point of $W^{n-1}$ corresponding to $\kappa$ lies
on the axis of revolution for some but not all points of  $W^{n-1}$, then $M^n$ is not proper Dupin.

Another problem with these constructions is that they may not yield
an immersed hypersurface in ${\bf R}^{n+1}$. In the tube construction, if the
radius of the tube is the reciprocal of one of the
principal curvatures of $W^{n-1}$ at some point, then the constructed object has a singularity.  For the
surface of revolution construction, a singularity occurs
if the prof\/ile submanifold $W^{n-1}$ intersects the axis of revolution. These problems can be resolved by working
in the context of Lie sphere geometry
(see \cite[pp.~127--141]{Cec1}).

As noted above, these constructions can be generalized to the setting of Lie sphere geometry by considering
Legendre lifts of hypersurfaces in Euclidean space.  In that context, a proper Dupin submanifold
$\lambda : M^{n-1} \rightarrow \Lambda ^{2n-1}$ is said to be {\em reducible} if it is
is locally Lie equivalent to the Legendre lift of a hypersurface in ${\bf R}^n$ obtained by one of Pinkall's constructions.

Pinkall \cite{P3} found a useful characterization of reducibility in the context of Lie sphere
geometry when he proved that a proper Dupin submanifold
$\lambda : M^{n-1} \rightarrow \Lambda ^{2n-1}$ is reducible if and only if
the image of one its curvature sphere maps $K$ lies in a linear subspace of codimension two in ${\bf P}^{n+2}$.
Much of the recent research in the f\/ield has focused on the classif\/ication of irreducible proper Dupin hypersurfaces.

The following example (see \cite[pp.~132--133]{Cec1} for more detail) shows
that Pinkall's constructions can produce a proper Dupin submanifold with the same number of distinct
curvature spheres as the original proper Dupin submanifold, rather than increasing the number of
distinct curvature spheres by one.

\looseness=-1
Let $V^2$ be a Veronese surface embedded in the unit sphere $S^4 \subset {\bf R}^5$.
Let $N^3$ be a tube of radius $\varepsilon$, for $0 < \varepsilon < \pi/3$, over $V^2$ in $S^4$.  Then $N^3$
is an isoparametric hypersurface $S^4$ with three distinct principal curvatures, and $N^3$ is an irreducible proper Dupin
hypersurface, since Cecil, Chi and Jensen \cite{CCJ2} proved that a compact, connected proper Dupin hypersurface
$M^{n-1} \subset S^n$ with $g \geq 3$ principal curvatures is irreducible.
$N^3$ is not the result of the tube construction, because
the Veronese surface is substantial (does not lie in a hyperplane) in ${\bf R}^5$.
Now embed ${\bf R}^5$ as a hyperplane through the origin in
${\bf R}^6$, and let $e_6$ be a unit normal vector to ${\bf R}^5$ in ${\bf R}^6$.
The Veronese surface $V^2$ is a codimension three submanifold of the
unit sphere $S^5 \subset {\bf R}^6$.  A tube $M^4$ over $V^2$ in $S^5$ is a reducible Dupin hypersurface, because the
family of curvatures spheres coming from the tube construction lies in a linear space of codimension two in ${\bf P}^7$,
since these curvature spheres all have the same radius and their centers all lie in the space ${\bf R}^5$.
This tube~$M^4$ is Dupin but not proper Dupin.
At points of the tube~$M^4$ coming from points of the form $(x,\pm e_6)$ in the unit normal
bundle $B^4$ of $V^2$ in~$S^5$,
there are two distinct principal curvatures, both of multiplicity
two.  At the other points of~$M^4$, there are three distinct principal curvatures, one of multiplicity two,
and the other two of multiplicity one.  Thus, the open dense subset $U$ of~$M^4$ on which there are
three principal curvatures is a reducible proper Dupin hypersurface, but~$M^4$ itself is not proper Dupin.
The number $g=3$ of distinct principal curvatures (or curvature spheres)
on $U$ is the same as the number of distinct curvature
spheres of the Legendre lift $\lambda$ of the original submanifold $V^2 \subset S^4$, since the point sphere
map of $\lambda$ is a~curvature sphere map, due to the fact that $V^2$ has codimension greater than one in $S^4$.
The other two curvature spheres of $\lambda$ are determined by the two principal curvatures of the surface~$V^2$.

\looseness=-1
This example is important for illustrating some of the subtleties involved in
studying the concept of reducibility.  For example,
Dajczer, Florit and Tojeiro \cite{DFT} studied reducibility in the context of Euclidean
submanifold theory. They formulated a concept of weak reducibility for proper Dupin submanifolds that have a
f\/lat normal bundle. In particular,
they def\/ined a~proper Dupin hypersurface
$f: M^{n-1} \rightarrow {\bf R}^n$ (or $S^n$)
to be {\em weakly reducible} if, for some principal curvature
$\kappa_i$ with corresponding principal space $T_i$, the orthogonal complement $T_i^{\perp}$ is integrable.
Dajczer, Florit and Tojeiro showed that if a proper Dupin hypersurface
$f: M^{n-1} \rightarrow {\bf R}^n$ is Lie equivalent to
a proper Dupin hypersurface with $g+1$ distinct principal curvatures that is obtained
from a proper Dupin hypersurface with $g$ distinct principal curvatures by one of the standard constructions,
then $f$ is weakly reducible.  Thus, reducible
implies weakly reducible for such hypersurfaces.
However, one can show that the open set $U$ with three principal curvatures in the example above
is reducible but not weakly reducible, because none
of the orthogonal complements of the principal spaces is integrable.
Note that $U$ is not constructed from a proper
Dupin submanifold with two curvature spheres, but rather one from one with three curvature spheres.

Cecil and Jensen \cite{CJ2,CJ3} def\/ined
a proper Dupin hypersurface $M^{n-1}$ in ${\bf R}^n$
to be {\em locally irreducible}
if $M^{n-1}$ does not contain a reducible open subset.
A locally irreducible proper Dupin hypersurface is obviously irreducible, and
using the analyticity of proper Dupin hypersurfaces, Cecil, Chi and Jensen \cite{CCJ2} proved conversely that
an irreducible proper Dupin hypersurface is locally irreducible. Thus, the two concepts are equivalent.

Using his constructions listed above, Pinkall \cite{P3} (see also \cite[p.~126]{Cec1})
demonstrated that proper Dupin hypersurfaces are plentiful.
Specif\/ically, he showed that given positive integers $m_1, \ldots , m_g$ with
$m_1 + \cdots + m_g = n-1$, there exists a proper
Dupin hypersurface $M^{n-1}$ in~${\bf R}^n$ with~$g$
distinct principal curvatures having respective
multiplicities $m_1, \ldots , m_g$.  The proper Dupin hypersurfaces that Pinkall constructs in his proof are all reducible,
and in general, they cannot be completed to a compact proper Dupin hypersurface, because of the dif\/f\/iculties
discussed above.  In fact, compact proper Dupin hypersurfaces are much more rare, as we will describe below.

Recall from Section~\ref{sec:1} that an immersion $\phi:V \rightarrow {\bf R}^n$ of a compact, connected manifold~$V$
into~${\bf R}^n$
is taut if every non-degenerate Euclidean distance function $L_p$, $p \in {\bf R}^n$, has the minimum
number of critical points on $V$ required by the Morse inequalities using ${\bf Z}_2$-homology.
Tautness can also be studied for submanifolds
of $S^n$ using spherical distance functions.

A taut submanifold of $S^n$
(or ${\bf R}^n$) must be Dupin, although not necessarily proper Dupin.  This was f\/irst shown for compact surfaces in
${\bf R}^3$ by Banchof\/f \cite{Ban1} (see also \cite{CecJDG} for the non-compact case),
then by Miyaoka \cite{Mi8} for hypersurfaces, and independently by Pinkall \cite{P4}
for submanifolds of arbitrary codimension.

Conversely, Thorbergsson \cite{Th1} proved that a compact, connected proper Dupin hypersurface $M^{n-1} \subset S^n$
(or ${\bf R}^n$) with $g$ distinct principal curvatures is taut.
Thorbergsson used the principal foliations on $M^{n-1}$ to construct concrete ${\bf Z}_2$-cycles in $M^{n-1}$ which
show that all critical points of non-degenerate distance functions $L_p$ are of linking type
(see Morse and Cairns \cite[p.~258]{MC}). Therefore
such distance functions have the minimum number $\beta (M^{n-1}) = 2g$ critical points, where $\beta (M^{n-1})$ is
the sum of the ${\bf Z}_2$-Betti numbers of $M^{n-1}$.
Pinkall \cite{P4} later
extended Thorbergsson's result to compact submanifolds of higher
codimension for which the number of distinct principal curvatures is constant on the unit
normal bundle.

Note that unlike the case for isoparametric hypersurfaces, the focal submanifolds of a
compact, connected proper Dupin hypersurface $M^{n-1}$ need not be taut.  For example, for a ring cyclide of
Dupin $M^2 \subset {\bf R}^3$ obtained by inverting a torus of revolution in a sphere, one of the focal
submanifolds is an ellipse, which is not a taut embedding of $S^1$.

\looseness=1
Ozawa \cite{Oz} proved that if $V$ is a taut compact submanifold of ${\bf R}^n$ (or $S^n$),
then every connected component $S$ of
a critical set of a distance function $L_p$ on $V$ is a smooth submanifold of $V$, which is non-degenerate as a
critical submanifold in the sense of Morse--Bott critical point theory (see~\cite{BS}),
and $S$ itself is taut in ${\bf R}^n$.  Using Ozawa's result, one can prove (see \cite[p.~154]{CecMSRI}) that if
$M^{n-1} \subset {\bf R}^n$ is a taut compact, connected hypersurface, then given any
principal space~$T_{\mu}$ at any point $x \in M^{n-1}$, there is a curvature surface $S$ through $x$ whose tangent space
is equal to~$T_{\mu}$, and $\mu$ is constant along $S$.  Thus $M^{n-1}$ is Dupin in this strong sense of having a curvature
surface tangent to every principal space (this is not assumed as part of Condition (a) in the def\/inition of a Dupin
hypersurface). An important open question is whether the converse of this theorem is true. That is, if
$M^{n-1} \subset {\bf R}^n$ is
a compact, connected non-proper Dupin hypersurface with the property that given any
principal space~$T_{\mu}$ at any point $x \in M^{n-1}$, there is a curvature surface $S$ through $x$ whose tangent space
is equal to $T_{\mu}$, must $M^{n-1}$ be taut?  Thorbergsson's
proof that a compact proper Dupin hypersurface must be taut
relies on the fact that all the curvature surfaces are spheres.  In the non-proper Dupin
case, Ozawa's work implies that there are some curvature surfaces that are not spheres.

Terng def\/ined a Dupin submanifold $V$ of arbitrary codimension
to have {\em constant multiplici\-ties} if
the multiplicities of the principal curvatures of any parallel normal
f\/ield $\xi(t)$ along any piecewise smooth curve on $V$ are constant.  Terng \cite{Te3}, \cite[p.~467]{Te2} then
proved that a compact, connected Dupin submanifold with constant multiplicities is taut.

Using tautness, Thorbergsson \cite{Th1} proved that a compact, connected proper Dupin hypersurface $M^{n-1} \subset S^n$
separates $S^n$ into two disk bundles over the f\/irst focal submanifolds on either side of $M^{n-1}$, as in the
case of isoparametric hypersurfaces. Thus, M\"{u}nzner's restriction that the number of distinct principal curvatures must be
 1, 2, 3, 4  or 6 holds for compact, connected proper Dupin hypersurfaces $M^{n-1}\subset S^n$ also.
Furthermore, the
restrictions on the possible multiplicities of the principal curvatures due to Stolz~\cite{Stolz} in the case
$g=4$ and Grove and Halperin~\cite{GH} in the case $g=6$ still hold.  These results have led to a pursuit of classif\/ication
results for compact proper Dupin hypersurfaces based on the number $g$ of distinct principal curvatures.

In the case $g=1$, it is well known that a compact, connected proper Dupin hypersurface~$M^{n-1} \subset S^n$ must be
a great or small hypersphere.  In the case $g=2$, Cecil and Ryan~\cite{CRMA} (see also~\cite{CRCan} or
\cite[p.~168]{CR}) showed that $M^{n-1}$ must be
M\"{o}bius equivalent to an isoparametric hypersurface, i.e., a standard product of spheres
$S^k(r) \times S^{n-k-1}(s) \subset S^n$, $r^2+s^2=1$. In the case $g=3$, Miyaoka \cite{Mi1} proved that
a compact, connected proper Dupin hypersurface $M^{n-1} \subset S^n$ must be Lie equivalent to
an isoparametric hypersurface, although not necessarily M\"{o}bius equivalent.

These results together with Thorbergsson's restriction on the number $g$ of distinct principal curvatures led to the
widely held conjecture \cite[p.~184]{CR} that every compact, connected
proper Dupin hypersurface $M^{n-1} \subset S^n$ is Lie
equivalent to an isoparametric hypersurface in a sphere.  However, in 1988
Pinkall and Thorbergsson~\cite{PT1}, and Miyaoka and Ozawa~\cite{MO} gave two dif\/ferent methods for producing
counterexamples to this conjecture with four principal curvatures.  The method of Miyaoka and Ozawa also yields
counterexamples to the conjecture in the case of six principal curvatures.

\looseness=1
Both of these constructions involve a consideration of the Lie curvatures of a proper Dupin hypersurface.
If $M^{n-1} \subset S^n$ is a proper Dupin hypersurface with $g \geq 4$ distinct principal curvatures, then
the cross-ratios of the principal curvatures taken four at a time are called the {\em Lie curvatures} of $M^{n-1}$.
These Lie curvatures were f\/irst studied by Miyaoka \cite{Mi2} who proved that they
are invariant under Lie sphere transformations.  This is actually
quite easy to see in the context of projective geometry, since each Lie curvature is the cross-ratio
of four points (corresponding to curvature spheres) on a projective line in ${\bf P}^{n+2}$.  A Lie sphere
transformation~$\beta$ is a projective transformation and it maps curvature spheres of a Legendre
submanifold $\lambda$ to curvature spheres of the Legendre submanifold $\beta \lambda$, so it preserves
the cross-ratios of the curvature spheres, and therefore it preserves the cross-ratios of the principal curvatures
(see also \cite[p.~75]{Cec1}).

The counterexamples of Pinkall and Thorbergsson \cite{PT1} (see also \cite[pp.~112--117]{Cec1})
to the conjecture are obtained by modifying the isoparametric hypersurfaces
of FKM-type discussed in the previous section. Recall that the Clif\/ford--Stiefel manifold
$V_2 (C_{m-1})$ is a submanifold of $S^{2l-1} \subset {\bf R}^{2l} = {\bf R}^l \times {\bf R}^l$
of codimension $m+1$.
Given positive real numbers $\alpha$ and $\beta$ with $\alpha^2 + \beta^2 = 1$, where
$\alpha \neq 1/\sqrt{2}$, $\beta \neq 1/\sqrt{2}$, Pinkall and Thorbergsson def\/ine a linear map
$T_{\alpha,\beta} : {\bf R}^{2l} \rightarrow {\bf R}^{2l}$, by
$T_{\alpha,\beta} (u,v) = \sqrt{2}\,  (\alpha u, \beta v)$.
Then for $(u,v) \in V_2 (C_{m-1})$, we have
\begin{gather*}
|T_{\alpha,\beta} (u,v)|^2 = 2 \big(\alpha^2 (u \cdot u) + \beta^2 (v \cdot v)\big) = 2 \left(\frac{\alpha^2}{2} + \frac{\beta^2}{2}\right)=1,
\end{gather*}
and thus the image $T_{\alpha,\beta}V_2 (C_{m-1})$ is a submanifold of $S^{2l-1}$ of
codimension $m+1$ also.  Pinkall and Thorbergsson prove that a tube over
$T_{\alpha,\beta}V_2 (C_{m-1})$ in $S^{2l-1}$ is a
compact, connected
proper Dupin hypersurface with four principal curvatures that does not have constant Lie curvature, and therefore it
is not Lie equivalent to an isoparametric hypersurface.

\looseness=1
The construction of counterexamples to the conjecture due to Miyaoka and Ozawa \cite{MO} (see also
\cite[pp.~117--123]{Cec1}) is based on the Hopf f\/ibration
$h:S^7 \rightarrow S^4$.
Miyaoka and Ozawa f\/irst show that if $M$ is a taut
compact submanifold of~$S^4$, then $h^{-1}(M)$ is a taut compact submanifold of
$S^7$.  Using this, they next show that if $M$ is a proper Dupin hypersurface in~$S^4$ with $g$ distinct principal curvatures, then
$h^{-1}(M)$ is a proper Dupin hypersurface in~$S^7$ with $2g$ principal
curvatures.  To complete the argument, they show that if a hypersurface $M \subset S^4$ is proper Dupin
but not isoparametric, then the Lie curvatures of
$h^{-1}(M)$ are not constant, and therefore $h^{-1}(M)$ is not Lie
equivalent to an isoparametric hypersurface in $S^7$.  For $g=2$ or 3, respectively, this gives
a compact proper Dupin hypersurface with 4 or 6
principal curvatures, respectively, in $S^7$, which is not Lie equivalent to
an isoparametric hypersurface.

These counterexamples are both based on the fact that the constructed proper Dupin hypersurface does not have
constant Lie curvatures.  This leads to a revision of the conjecture \cite[p.~52]{CCJ4} which states
that every compact, connected proper Dupin hypersurface $M^{n-1} \subset S^n$ with $g=4$ or 6 principal curvatures
and constant Lie curvatures is Lie equivalent to an isoparametric hypersurface in a sphere.  This revised conjecture
is still an important open problem, although it has has been shown to be true in some cases, as we now describe.

Miyaoka \cite{Mi2,Mi3} began by showing that if some additional assumptions are made regarding the intersections
of the leaves of the various principal foliations, then this revised conjecture is true in both cases
$g=4$ and $g=6$.
Thus far, however, it has not been shown that Miyaoka's additional assumptions are satisf\/ied in general. In the
case $g=6$, this work of Miyaoka is the only known progress towards proving the revised conjecture.

If $M^{n-1} \subset S^n$ is a compact, connected proper
Dupin hypersurface with four principal curvatures having multiplicities $m_1$, $m_2$, $m_3$, $m_4$, then
the multiplicities must satisfy
$m_1 = m_2$, $m_3 = m_4$, when the principal curvatures are appropriately ordered, by the work of Thorbergsson~\cite{Th1}
and M\"{u}nzner \cite{Mu,Mu2}.
Cecil, Chi and Jensen \cite{CCJ2} have recently shown
that if $M^{n-1} \subset S^n$ is a compact, connected proper
Dupin hypersurface with four principal curvatures having multiplicities
$m_1 = m_2 \geq 1$, $m_3 = m_4 = 1$, and constant Lie curvature,
then $M^{n-1}$ is Lie equivalent to an isoparametric hypersurface.
This result actually follows from a
local classif\/ication of irreducible proper Dupin hypersurfaces with four principal curvatures~\cite{CCJ2}, which we will
discuss below. Thus, the full revised conjecture in the case $g=4$ would be proven if one could remove the assumption
$m_3 = m_4 = 1$, but so far this has not been done.

In his early work on the subject, Pinkall \cite{P1,P2,P3}, proved some important local classif\/ication
results concerning proper Dupin hypersurfaces, and this approach has been extended to more general settings by using
the notion of irreducibility.
Although it has not been shown that irreducibility places any restriction on the number $g$ of distinct
principal curvatures other than $g \geq 3$, most of the known results for irreducible proper Dupin
hypersurfaces have been obtained in the cases $g = 3$ or~4.

As with compact proper Dupin hypersursurfaces, we will discuss the known local classif\/ication results
for proper Dupin hypersurfaces of $S^n$ based on the number $g$ of distinct principal curvatures.  It is well known,
that a connected proper Dupin hypersurface in $S^n$ with $g=1$ (totally umbilic) distinct principal curvature
must be an open subset of a hypersphere.

In the case $g=2$, Pinkall \cite{P2} obtained a complete local classif\/ication which we now describe in detail.
As a generalization of the well-known cyclides of Dupin in ${\bf R}^3$, Pinkall def\/ined a~{\em cyclide of Dupin of characteristic} $(p,q)$ to be a proper Dupin hypersurface in $S^n$ (or ${\bf R}^n$)
with two distinct principal curvatures of respective multiplicities $p$ and $q$.  An example of a
cyclide of Dupin of characteristic $(p,q)$ is the standard product of spheres, $S^p (r) \times S^q (s) \subset S^n$,
$r^2 + s^2 = 1$, where $n = p+q+1$, as in equation (\ref{product}).
This is an isoparametric hypersurface in $S^n$ with two principal curvatures, as discussed earlier.
If one varies the value of $r$ in equation (\ref{product}) between~0 and~1, one obtains an isoparametric family
of hypersurfaces with two principal curvatures.  The hypersurfaces in this family are Lie equivalent
by parallel transformation, but they are not M\"{o}bius equivalent for dif\/ferent values of $r$.

Pinkall proved that every connected cyclide of Dupin is contained in a
unique compact, connected cyclide, and any two cyclides of the same characteristic $(p,q)$ are locally
Lie equivalent. Thus, every cyclide of Dupin of characteristic $(p,q)$ is locally Lie equivalent to a standard
product of spheres, as in equation (\ref{product}).

Using Pinkall's Lie geometric classif\/ication one can derive the following
M\"{o}bius geometric classif\/ication of the cyclides of Dupin \cite{CecL} (see also \cite[p.~151]{Cec1}).
Every connected cyclide $M^{n-1}$ of characteristic
$(p,q)$ in ${\bf R}^n$ is M\"{o}bius equivalent to an open
subset of a hypersurface of revolution obtained by revolving a
$q$-sphere $S^q \subset {\bf R}^{q+1} \subset {\bf R}^n$ about an axis
of revolution ${\bf R}^q \subset {\bf R}^{q+1}$ or a $p$-sphere $S^p
\subset {\bf R}^{p+1} \subset {\bf R}^n$ about an axis
${\bf R}^p \subset {\bf R}^{p+1}$. Further, two such hypersurfaces of revolution are M\"{o}bius equivalent if and only if
they have the same value of $\rho = r/a$, where $r$ is the
radius of the prof\/ile sphere $S^q$ and $a>0$ is the distance from the center of $S^q$
to the axis of revolution.
In this theorem, the prof\/ile sphere is allowed to intersect the axis
of revolution, which results in Euclidean singularities.  However,
in the context of Lie sphere geometry, the corresponding Legendre map is an immersion.

The classical cyclides of Dupin in ${\bf R}^3$ were f\/irst studied by
Dupin\index{Dupin} \cite{D} in 1822 and later by many prominent nineteenth century mathematicians,
including Liouville \cite{Lio}, Cayley \cite{Cay}, and Maxwell \cite{Max}, whose paper contains stereoscopic
f\/igures of the various types of cyclides. (See Lilienthal \cite{Lil} for an account
of the nineteenth century work on the cyclides.)
More recent descriptions of the classical cyclides are contained in the books of
Fladt and Baur \cite[pp.~354--379]{FB}, Cecil and Ryan
\cite[pp.~151--166]{CR}, and \cite[pp.~148--159]{Cec1}.

The classical cyclides are the only surfaces in ${\bf R}^3$ with two principal curvatures at each point
such that all lines of curvature in both families are circles or straight lines. This is just the proper Dupin condition,
of course.
Using exterior dif\/ferential systems,
Ivey \cite{Ivey} showed that any surface in ${\bf R}^3$ containing
two orthogonal families of circles is a cyclide of Dupin.
The classical cyclides have also been
used in computer aided geometric design of surfaces.  See, for example, the papers of
Degen \cite{Degen}, Pratt \cite{Pr1,Pr2},
Srinivas and Dutta \cite{SD1,SD2,SD3,SD4}, and
Schrott and Odehnal \cite{Schrott}.

Pinkall began the study of proper Dupin hypersurfaces with three distinct principal curvatures in his dissertation \cite{P1},
published as a paper \cite{P3} (see also \cite{CC2} or \cite[pp.~168--190]{Cec1}).  Pinkall found
a complete local classif\/ication up to Lie equivalence for Dupin hypersurfaces
with three principal curvatures in ${\bf R}^4$.  He proved that any two irreducible
proper Dupin hypersurfaces with $g=3$ in ${\bf R}^4$ are locally
Lie equivalent, each being Lie equivalent to an open subset
of Cartan's isoparametric hypersurface in $S^4$.
For reducible proper Dupin hypersurfaces with $g=3$ in~${\bf R}^4$, Pinkall showed that there is
a 1-parameter family of Lie equivalence classes.

Niebergall \cite{N1} next proved that every connected proper Dupin hypersurface
in ${\bf R}^5$ with three principal curvatures is reducible.  Then
Cecil and Jensen \cite{CJ2} proved that if $M^{n-1}$ is an irreducible, connected proper Dupin hypersurface in~$S^n$ with
three distinct principal curvatures of multiplicities
$m_1$, $m_2$, $m_3$, then $m_1 = m_2 = m_3 = m$, and $M^{n-1}$ is Lie
equivalent to an isoparametric hypersurface in $S^n$.  It then follows from Cartan's classif\/ication of
isoparametric hypersurfaces with $g=3$ mentioned in Section~\ref{sec:1} that $m=1,2,4$ or 8.
Note that in the original paper~\cite{CJ2}, this result was proven under the assumption that
$M^{n-1}$ is locally irreducible.  However, as noted above, local irreducibility has now been shown to be equivalent to
irreducibility.

The proof of this result of Cecil and Jensen is accomplished by using Cartan's method of moving frames in the
context of Lie sphere geometry.  A key tool in this context is the result (see~\cite{CecKod} or \cite[p.~77]{Cec1})
that a Legendre submanifold $\lambda : M^{n-1} \rightarrow \Lambda ^{2n-1}$
with $g$ distinct curvature spheres
$K_1,\ldots,K_g$ at each point is
Lie equivalent to the Legendre lift of an
isoparametric hypersurface in $S^n$ if and only if there
exist $g$ points $P_1,\ldots,P_g$ on a timelike line in
$P^{n+2}$ such that $\langle K_i,P_i\rangle = 0$, for $1 \leq i \leq g$.

An open problem is the classif\/ication of reducible Dupin hypersurfaces
with three principal curvatures up to Lie equivalence.
As noted above, Pinkall \cite{P3} found such a classif\/ication in the
case of $M^3 \subset {\bf R}^4$. It may be
possible to generalize Pinkall's result to higher
dimensions using the approach of \cite{CJ2}.

Cecil, Chi and Jensen \cite{CCJ2} proved that a compact, connected proper Dupin hypersurface $M^{n-1} \subset S^n$
with $g \geq 3$ principal curvatures is irreducible.  Thus, the examples above
of Pinkall and Thorbergsson \cite{PT1} and
Miyaoka and Ozawa \cite{MO} of compact proper Dupin hypersurfaces with $g=4$ and non-constant
Lie curvature are irreducible proper Dupin
hypersurfaces with $g=4$ that are not Lie equivalent to an isoparametric hypersurface.  However, it is still possible that
every irreducible Dupin hypersurface with $g=4$ and constant Lie curvature
is Lie equivalent to an isoparametric hypersurface, and this has been shown under additional assumptions,
as we now describe.

The f\/irst result in this direction is due to Niebergall \cite{N2} who showed that a connected irreducible proper Dupin
hypersurface $M^4$ in $S^5$ with four principal curvatures and constant Lie curvature is Lie equivalent to an
isoparametric hypersurface under an additional assumption that in an appropriate moving frame, the covariant
derivatives of certain naturally def\/ined functions are zero.  Cecil and Jensen \cite{CJ3} then showed that Niebergall's
additional assumptions are unnecessary, because these functions are always constant in the appropriate
Lie frame.  Thus, they proved that every connected irreducible proper Dupin hypersurface
$M^4$ in $S^5$ with four principal curvatures and constant Lie curvature is Lie equivalent to an
isoparametric hypersurface.

This result was then generalized to higher dimensions by Cecil, Chi and Jensen \cite{CCJ2} who
proved that if $M^{n-1}$ is an irreducible, connected proper Dupin hypersurface in $S^n$ with
four principal curvatures having multiplicities $m_1=m_2\geq 1$ and $m_3=m_4=1$ and constant Lie curvature $-1$,
then $M^{n-1}$ is Lie equivalent to an isoparametric hypersurface.
Note that M\"{u}nzner \cite{Mu,Mu2} had shown earlier that if $M^{n-1} \subset S^n$ is an isoparametric
hypersurface with four principal curvatures, then the multiplicities of the principal curvatures must satisfy
$m_1=m_2$ and $m_3=m_4$, and the Lie curvature must equal $-1$, if the principal curvatures are appropriately ordered.

These results lead to the following local conjecture \cite[p.~53]{CCJ4}:
if~$M^{n-1}$ is an irreducible, connected, proper Dupin hypersurface in~$S^n$
with four principal curvatures having multiplicities $m_1$, $m_2$, $m_3$, $m_4$, and constant Lie curvature $c$, then
$m_1=m_2$, $m_3=m_4$, $c=-1$, and $M^{n-1}$ is Lie equivalent to an isoparametric hypersurface.
Note that it has not yet been shown that irreducibility implies that $m_1=m_2$ and $m_3=m_4$, nor that $c=-1$.

\looseness=1
We remark that the hypothesis of irreducibility is def\/initely needed in the conjecture because one can construct reducible
non-compact proper Dupin hypersurfaces
with $g=4$ and constant Lie curvature $c$, for every negative value of $c$, if the principal curvatures
are appropriately ordered (see \cite{CecKod} or \cite[pp.~80--82]{Cec1}).
This construction yields examples where the multiplicities satisfy
$m_1=m_2$ and $m_3=m_4$, and also examples where the multiplicities do not satisfy this restriction.
These examples are all obtained as open subsets of a tube in
$S^n$ over an isoparametric hypersurface with three principal
curvatures $V^{k-1} \subset S^k \subset S^n$, and they
cannot be completed to be compact proper Dupin hypersurfaces.
They are also reducible as Dupin hypersurfaces.

Using a dif\/ferent approach based on the theory of higher-dimensional Laplace invariants \cite{KT},
Riveros and Tenenblat \cite{RT,RT2} gave a local classif\/ication of
proper Dupin hypersurfaces $M^4$ in ${\bf R}^5$ with four distinct principal curvatures which
are parametrized by lines of curvature.

C.-P. Wang \cite{Wc,Wc2} studied the M\"{o}bius geometry of submanifolds in $S^n$
in a series of papers. Using Cartan's method of moving frames, Wang found a complete set of M\"{o}bius invariants
for surfaces in ${\bf R}^3$ without umbilic points \cite{Wc} and for hypersurfaces in ${\bf R}^4$ with three
principal curvatures at each point~\cite{Wc1}.  Then
in~\cite{Wc2}, Wang def\/ined a M\"{o}bius invariant metric $g$ and second fundamental form $B$
for submanifolds in $S^n$.  Wang then proved that for hypersurfaces in $S^n$ with $n \geq4$,
the pair $(g,B)$ forms a complete
M\"{o}bius invariant system which determines the hypersurface up to M\"{o}bius transformations.

In a related result, Riveros, Rodrigues and Tenenblat \cite{RRT} proved that a proper Dupin hypersurface
$M^n \subset {\bf R}^{n+1}$, $n \geq 4$, with $n$ distinct principal curvatures and constant M\"{o}bius curvatures
cannot be parametrized by lines of curvature. They also showed that up to M\"{o}bius transformations, there is a
unique proper Dupin hypersurface $M^3 \subset {\bf R}^4$ with three principal curvatures and constant
M\"{o}bius curvature that is parametrized by lines of curvature. This $M^3$ is a cone in ${\bf R}^4$ over a
standard f\/lat torus in the unit sphere $S^3 \subset {\bf R}^4$.

In \cite{LLWZ},
H. Li, Lui, Wang and Zhao introduced the concept of a M\"{o}bius isoparametric hypersurface in a sphere $S^n$.
They showed that a (Euclidean) isoparametric hypersurface is automatically M\"{o}bius isoparametric, whereas a
M\"{o}bius isoparametric hypersurface must be proper Dupin. Later Rodrigues and Tenenblat \cite{RoT}
showed that if $M \subset S^n$ is a hypersurface with a~constant number $g$
of distinct principal curvatures at each point, where $g \geq 3$, then
$M$ is M\"{o}bius isoparametric
if and only if $M$ is Dupin with constant M\"{o}bius curvatures.

Recently signif\/icant progress has been made in the classif\/ication of M\"{o}bius isoparametric hypersurfaces. First,
H. Li, Lui, Wang and Zhao \cite{LLWZ} showed that a connected M\"{o}bius
isoparametric hypersurface in $S^n$ with two distinct
principal curvatures is M\"{o}bius equivalent to an open subset
of one of the following three types of hypersurfaces in $S^n$:

\begin{enumerate}\itemsep=0pt
\item[(i)] a standard product of spheres $S^k(r) \times S^{n-k-1}(s) \subset S^n$, $r^2+s^2=1$,
\item[(ii)] the image under inverse stereographic projection from ${\bf R}^n \rightarrow S^n - \{P\}$ of a standard cylinder
$S^k(1) \times {\bf R}^{n-k-1} \subset {\bf R}^n$,
\item[(iii)] the image under hyperbolic stereographic projection from $H^n \rightarrow S^n$ of a standard product
$S^k(r) \times H^{n-k-1}(\sqrt{1+r^2}) \subset H^n$.
\end{enumerate}

\noindent
Later Hu and H.~Li \cite{HL-H} classif\/ied M\"{o}bius isoparametric hypersurfaces in $S^4$, and Hu, H. Li and Wang \cite{HLW}
classif\/ied M\"{o}bius isoparametric hypersurfaces in $S^5$.  Then Hu and D.~Li \cite{HL-D}
studied M\"{o}bius isoparametric hypersurfaces with three distinct principal curvatures in~$S^n$ and found a~complete classif\/ication of such hypersurfaces in $S^6$.

\subsection*{Acknowledgements}

This material is based upon work supported by the National Science Foundation under Grant No.~0405529.
The author is grateful for several helpful comments in the reports of the referees.

\pdfbookmark[1]{References}{ref}
\LastPageEnding


\begin{thebibliography}{99}

\footnotesize\itemsep=0pt
\bibitem{Ab} Abresch U., Isoparametric hypersurfaces with
four or six distinct principal curvatures, {\em Math. Ann.} {\bf 264} (1983),
283--302.

\bibitem{Alex} Alexandrino M., Singular Riemannian foliations with sections,
{\em Illinois J. Math.} {\bf 48} (2004), 1163--1182, \href{http://arxiv.org/abs/math.DG/0311454}{math.DG/0311454}.

\bibitem{ABS} Atiyah M.F., Bott R., Shapiro A., Clif\/ford modules, {\em Topology} {\bf 3} (1964), suppl.~1, 3--38.

\bibitem{Ban1} Banchof\/f T., The spherical two-piece property and
tight surfaces in spheres, {\em J. Differential Geom.} {\bf 4} (1970), 193--205.

\bibitem{Berndt} Berndt J., Real hypersurfaces with constant principal curvatures in complex space forms,
in Proceedings of the Tenth International Workshop on Dif\/ferential Geometry, Kyungpook Nat. Univ., Taegu, 2006, 1--12.

\bibitem{BCO} Berndt J., Console  S., Olmos C., Submanifolds and holonomy, {\it Chapman and Hall/CRC Research Notes
in Mathematics}, Vol.~434, Chapman and Hall/CRC, Boca Raton, Florida, 2003.

\bibitem{BS} Bott  R., Samelson H., Applications of the
theory of Morse to symmetric spaces, {\em Amer. J. Math.} {\bf 80} (1958),
964--1029.

\bibitem{Car1} Cartan \'E., Familles de surfaces
isoparam\'{e}triques dans les espaces \`{a} courbure constante,
{\em Annali di Mat.} {\bf 17} (1938), 177--191 (see also in {\em Oeuvres Compl\`{e}tes}, Partie III, Vol. 2, 1431--1445).

\bibitem{Car2} Cartan \'E., Sur des familles remarquables d'hypersurfaces isoparam\'{e}triques dans les espaces
sph\'{e}riques, {\em Math. Z.} {\bf 45} (1939), 335--367
(see also in {\em Oeuvres Compl\`{e}tes}, Partie III, Vol. 2, 1447--1479).

\bibitem{Car3} Cartan \'E., Sur quelque familles remarquables d'hypersurfaces,
in C.R. Congr\`{e}s Math. Li\`{e}ge, 1939, 30--41 (see also in {\em Oeuvres Compl\`{e}tes}, Partie III, Vol. 2, 1481--1492).

\bibitem{Car4} Cartan \'E., Sur des familles d'hypersurfaces isoparam\'{e}triques des espaces sph\'{e}riques \`{a} 5
et \`{a} 9 dimensions, {\em Revista Univ. Tucuman, Serie A}, {\bf 1} (1940), 5--22 (see also in {\em Oeuvres Compl\`{e}tes}, Partie III, Vol. 2, 1513--1530).

\bibitem{CS}Carter  S., \c{S}ent\"{u}rk Z., The space of immersions parallel to a given immersion,
{\em J. London Math. Soc. (2)} {\bf 50} (1994), 404--416.

\bibitem{CW} Carter  S., West A., Tight and taut immersions, {\em Proc. London Math. Soc. (3)} {\bf 25} (1972), 701--720.

\bibitem{CW3} Carter  S., West A., Totally focal embeddings, {\em J. Differential Geom.} {\bf 13} (1978), 251--261.

\bibitem{CW4} Carter  S., West A., Totally focal embeddings: special cases, {\em J. Differential Geom.} {\bf 16} (1981), 685--697.

\bibitem{CW6} Carter  S., West A., A characterisation of isoparametric hypersurfaces in spheres, {\em J. London Math. Soc. (2)} {\bf 26} (1982), 183--192.

\bibitem{CW1} Carter S., West A., Isoparametric
systems and transnormality, {\em Proc. London Math. Soc. (3)} {\bf 51} (1985), 520--542.

\bibitem{CW5} Carter  S., West A., Isoparametric and totally focal submanifolds, {\em Proc. London Math. Soc. (3)}
{\bf 60} (1990), 609--624.

\bibitem{Cay} Cayley A., On the cyclide, {\em Quart. J. of Pure and
Appl. Math.} {\bf 12} (1873), 148--165 (see also in {\em Collected Mathematical Papers}, Vol.~9, Cambridge U. Press, 1896, 64--78).

\bibitem{CecToh} Cecil T., A characterization of metric spheres in hyperbolic space by Morse theory,
{\em T\^{o}hoku Math. J. (2)} {\bf 26} (1974), 341--351.

\bibitem{CecJDG} Cecil T., Taut immersions of non-compact
surfaces into a Euclidean 3-space, {\em J. Differential Geom.} {\bf 11} (1976), 451--459.

\bibitem{CecGD} Cecil T., Reducible Dupin submanifolds, {\em Geom. Dedicata} {\bf 32} (1989), 281--300.

\bibitem{CecKod} Cecil T., On the Lie curvatures of Dupin hypersurfaces, {\em Kodai Math. J.} {\bf 13} (1990), 143--153.

\bibitem{CecL} Cecil T., Lie sphere geometry and Dupin submanifolds, in  Geometry and Topology of Submanifolds III
 (Leeds, 1990),
Editors L.~Verstraelen and A.~West, World Scientif\/ic, River Edge, NJ, 1991, 90--107.

\bibitem{CecMSRI} Cecil T., Taut and Dupin submanifolds, in Tight and Taut Submanifolds (Berkeley, CA, 1994), Editors T.~Cecil and
S.-S.~Chern, {\it Math. Sci. Res. Inst. Publ.}, Vol.~32, Cambridge Univ. Press, Cambridge, 1997, 135--180.

\bibitem{Cec1} Cecil T., Lie sphere geometry, with applications to submanifolds, 2nd ed.,
Universitext, Springer, New York, 2008.

\bibitem{CC1} Cecil T., Chern S.-S., Tautness and Lie sphere geometry, {\em Math. Ann.} {\bf 278} (1987), 381--399.

\bibitem{CC2} Cecil  T., Chern S.-S., Dupin submanifolds in Lie sphere geometry, in  Dif\/ferential Geometry and Topology (Tianjin, 1986--87), Editors B.~Jiang et al.,
{\it Lecture Notes in Math.}, Vol.~1369, Springer, Berlin~-- New York, 1989, 1--48.

\bibitem{CCJ1} Cecil T., Chi  Q.-S., Jensen G., Isoparametric hypersurfaces with four principal
curvatures, {\em Ann. of Math.~(2)} {\bf 166} (2007), 1--76. 

\bibitem{CCJ2} Cecil T.,  Chi Q.-S., Jensen G., Dupin hypersurfaces with four principal curvatures.~II,
{\em Geom. Dedicata} {\bf 128} (2007), 55--95.  

\bibitem{CCJ4} Cecil T., Chi Q.-S., Jensen G., Classif\/ications of Dupin hypersurfaces,
in  Pure and Applied Dif\/ferential Geometry, PADGE 2007, Editors F.~Dillen and I.~van de Woestyne,
Shaker Verlag, Aachen, 2007, 48--56.

\bibitem{CCJ3} Cecil T., Chi  Q.-S., Jensen G., On Kuiper's question whether taut submanifolds are algebraic,
{\em Pacific J. Math.} {\bf 234} (2008), 229--247.

\bibitem{CJ2} Cecil  T., Jensen G., Dupin hypersurfaces with three principal curvatures,
{\em Invent. Math.} {\bf 132} (1998), 121--178.

\bibitem{CJ3} Cecil  T.,  Jensen G., Dupin hypersurfaces with four principal curvatures,
{\em Geom. Dedicata} {\bf 79} (2000), 1--49.

\bibitem{CRMA} Cecil  T., Ryan P., Focal sets, taut embeddings and the cyclides of Dupin, {\em Math. Ann.}
{\bf 236} (1978), 177--190.

\bibitem{CRCan} Cecil  T., Ryan P., Conformal geometry and the cyclides of Dupin, {\em Canad. J. Math.}
{\bf 32} (1980), 767--782.

\bibitem{CRBer} Cecil  T., Ryan P., Tight spherical embeddings, in Global Dif\/ferential Geometry and Analysis (Berlin 1979), Editors D.~Ferus, W. K\"{u}hnel, U. Simon and B. Wegner,
{\it Lecture Notes in Math.}, Vol.~838, Springer, Berlin~-- New York, 1981, 94--104.

\bibitem{CR} Cecil  T., Ryan P., Tight and taut immersions of manifolds, {\it Research Notes in Math.}, Vol.~107, Pitman,
London, 1985.

\bibitem{Chern} Chern S.-S., An introduction to Dupin submanifolds, in Dif\/ferential Geometry: A Symposium in
Honour of Manfredo do Carmo (Rio de Janeiro, 1988), Editors H.B.~Lawson and K.~Tenenblat,
{\it Pitman Monographs Surveys Pure Appl. Math.}, Vol.~52, Longman Sci. Tech., Harlow, 1991, 95--102.

\bibitem{Chi} Chi Q.-S., Isoparametric hypersurfaces with four principal curvatures revisited, in
{\em Nagoya Math. J.}, to appear, \href{http://arxiv.org/abs/0803.1284}{arXiv:0803.1284}.

\bibitem{Christ} Christ U., Homogeneity of equifocal submanifolds, {\em J. Differential Geom.} {\bf 62} (2002), 1--15.

\bibitem{Dad} Dadok J., Polar coordinates induced by actions of compact Lie groups,
{\em Trans. Amer. Math. Soc.} {\bf 288} (1985), 125--137.

\bibitem{DFT} Dajczer M., Florit L., Tojeiro R., Reducibility of Dupin submanifolds,
{\em Illinois J. Math.} {\bf 49} (2005), 759--791.

\bibitem{Degen} Degen W., Generalized cyclides for use in CAGD, in  Computer-Aided Surface Geometry and Design:
the Mathematics of Surfaces IV
(Bath, 1990), Editor A.~Bowyer, {\it Inst. Math. Appl. Conf.
New Ser.}, Vol.~48, Oxford Univ. Press, New York, 1994, 349--363.

\bibitem{DN1} Dorfmeister  J., Neher E., An algebraic approach to isoparametric hypersurfaces.~I,
{\em T\^{o}hoku Math. J. (2)} {\bf 35} (1983), 187--224.\\
 Dorfmeister  J., Neher E., An algebraic approach to isoparametric hypersurfaces.~II,
 {\em T\^{o}hoku Math. J. (2)} {\bf 35} (1983), 225--247.

\bibitem{DN2} Dorfmeister  J., Neher E., Isoparametric triple systems of algebra type, {\em Osaka J. Math.}
{\bf 20} (1983), 145--175.

\bibitem{DN3} Dorfmeister  J., Neher E., Isoparametric triple systems of FKM-type.~I, {\em Abh. Math. Sem. Hamburg}
{\bf 53} (1983), 191--216.

\bibitem{DN4} Dorfmeister  J., Neher E., Isoparametric triple systems of FKM-type.~II,
{\em Manuscripta Math.} {\bf 43} (1983), 13--44.

\bibitem{DN} Dorfmeister  J., Neher E., Isoparametric hypersurfaces, case $g = 6$, $m = 1$,
{\em Comm. Algebra} {\bf 13} (1985), 2299--2368.

\bibitem{DN6} Dorfmeister  J., Neher E., Isoparametric triple systems with special $Z$-structure,
{\em Algebras Groups Geom.} {\bf 7} (1990), 21--94.

\bibitem{D} Dupin C., Applications de g\'{e}om\'{e}trie et de m\'{e}chanique, Bachelier, Paris, 1822.

\bibitem{ES} Eschenburg  J.-H., Schroeder V., Tits distance of Hadamard manifolds and isoparametric
hypersurfaces, {\em Geom. Dedicata} {\bf 40} (1991), 97--101.

\bibitem{Fang3} Fang F., Multiplicities of principal curvatures of isoparametric hypersurfaces, Preprint,
Max Planck Institut f\"{u}r Mathematik, Bonn, 1996.

\bibitem{Fang1} Fang F., On the topology of isoparametric hypersurfaces with four distinct principal curvatures,
{\em Proc. Amer. Math. Soc.} {\bf 127} (1999), 259--264.

\bibitem{Fang2} Fang F., Topology of Dupin hypersurfaces with six principal curvatures, {\em Math. Z.}
{\bf 231} (1999), 533--555.

\bibitem{Fera} Ferapontov E.V., Dupin hypersurfaces and integrable hamiltonian systems of hydrodynamic type, which
do not possess Riemann invariants, {\em Differential Geom. Appl.} {\bf 5} (1995), 121--152.

\bibitem{Fera1} Ferapontov E.V., Isoparametric hypersurfaces in spheres, integrable
nondiagonalizable systems of hydrodynamic type, and $N$-wave systems, {\em Differential Geom. Appl.}
{\bf 5} (1995), 335--369.

\bibitem{FKM} Ferus D., Karcher  H., M\"{u}nzner H.-F., Clif\/fordalgebren und neue isoparametrische
Hyperf\/l\"{a}chen, {\em Math. Z.} {\bf 177} (1981), 479--502.

\bibitem{FB} Fladt  K., Baur A., Analytische Geometrie spezieller Fl\"{a}chen und Raumkurven, Friedr.
Vieweg and Sohn, Braunschweig, 1975.

\bibitem{Geatti} Geatti  L., Gorodski C., Polar orthogonal representations of real reductive algebraic groups,
{\em J. Algebra}, to appear, \href{http://arxiv.org/abs/0801.0574}{arXiv:0801.0574}.

\bibitem{GH} Grove  K., Halperin S., Dupin hypersurfaces, group actions, and the double mapping cylinder,
{\em J. Differential Geom.} {\bf 26} (1987), 429--459.

\bibitem{Hahn} Hahn J., Isoparametric hypersurfaces in the pseudo-Riemannian space forms, {\em Math. Z.} {\bf 187}
(1984), 195--208.

\bibitem{Hahn2} Hahn J., Isotropy representations of semisimple symmetric spaces and homogeneous
hypersurfaces, {\em J. Math. Soc. Japan} {\bf 40} (1988), 271--288.

\bibitem{Har} Harle C., Isoparametric families of submanifolds, {\em Bol. Soc. Brasil Mat.} {\bf 13} (1982),  no.~2, 35--48.

\bibitem{HLiu} Heintze  E.,  Liu X., Homogeneity of inf\/inite dimensional isoparametric submanifolds,
{\em Ann. of Math. (2)} {\bf 149} (1999), 149--181, \href{http://arxiv.org/abs/math.DG/9901150}{math.DG/9901150}.

\bibitem{HOT} Heintze E., Olmos C., Thorbergsson G., Submanifolds with constant principal curvatures and
normal holonomy groups, {\em Internat. J. Math.} {\bf 2} (1991), 167--175.

\bibitem{HL} Hsiang  W.-Y., Lawson H. B., Jr., Minimal submanifolds of low cohomogeneity,
{\em J. Differential Geom.} {\bf 5} (1971), 1--38.

\bibitem{HPT} Hsiang W.-Y., Palais R., Terng C.-L., The topology of isoparametric submanifolds, {\em J. Differential Geom.} {\bf 27} (1988), 423--460.

\bibitem{HL-D} Hu Z., Li D., M\"{o}bius isoparametric hypersurfaces with three distinct principal curvatures,
{\em Pacific J. Math.} {\bf 232} (2007), 289--311.

\bibitem{HL-H} Hu  Z., Li H.-Z., Classif\/ication of M\"{o}bius isoparametric hypersurfaces in $S^4$,
{\em Nagoya Math. J.} {\bf 179} (2005), 147--162.

\bibitem{HLW} Hu Z., Li H.-Z., Wang C.-P., Classif\/ication of M\"{o}bius isoparametric hypersurfaces in $S^5$,
{\em Monatsh. Math.} {\bf 151} (2007), 201--222.

\bibitem{Im} Immervoll S., On the classif\/ication of isoparametric hypersurfaces with four distinct
principal curvatures in spheres, {\em Ann. Math.}, to appear.

\bibitem{Ivey} Ivey T., Surfaces with orthogonal families of circles, {\em Proc. Amer. Math. Soc.} {\bf 123}
(1995), 865--872.

\bibitem{KT} Kamran  N., Tenenblat K., Laplace transformation in higher dimensions,
{\em Duke Math. J.} {\bf 84} (1996), 237--266.

\bibitem{Lev} Levi-Civita T., Famiglie di superf\/icie
isoparametrische nell'ordinario spacio euclideo, {\em Atti.
Accad. naz. Lincei. Rend. Cl. Sci. Fis. Mat. Natur.} {\bf 26} (1937), 355--362.

\bibitem{LLWZ} Li H.-Z., Liu H.-L., Wang C.-P., Zhao G.-S., M\"{o}bius isoparametric hypersurfaces in $S^{n+1}$ with
two distinct principal curvatures, {\em Acta Math. Sin. (Engl. Ser.)} {\bf 18} (2002), 437--446.

\bibitem{LX} Li Z.-Q., Xie X.-H., Space-like isoparametric hypersurfaces in Lorentzian space forms,
{\it J. Nanchang Univ. Natur. Sci. Ed.} {\bf 28} (2004), 113--117
(English transl: {\em Front. Math. China} {\bf 1} (2006), 130--137).

\bibitem{Lil} Lilienthal R., Besondere Fl\"{a}chen, in  Encyklop\"{a}die
der Math. Wissenschaften, Vol.~III, Teubner, Leipzig, 1902--1927, 269--354.

\bibitem{Lio} Liouville J., Note au sujet de l'article pr\'{e}ced\'{e}nt, {J. de Math. Pure et Appl. (1)} {\bf 12}
(1847), 265--290.

\bibitem{LT} Lytchak  A., Thorbergsson G., Variationally complete actions on nonnegatively curved manifolds,
{\em Illinois J. Math.} {\bf 51} (2007), 605--615.

\bibitem{Mag} Magid M., Lorentzian isoparametric hypersurfaces, {\em Pacific J. Math.} {\bf 118} (1985), 165--197.

\bibitem{Max} Maxwell J.C., On the cyclide, {\em Quart. J. of Pure and Appl. Math.} {\bf 34} (1867), 111--126
(see also in {\em Scientific papers}, Vol.~2, Cambridge U. Press, 1890, 144--159).

\bibitem{Mil} Milnor J., Morse theory, {\em Ann. Math. Stud.}, Vol.~51, Princeton U. Press, Princeton, NJ, 1963.

\bibitem{Mi1} Miyaoka R., Compact Dupin hypersurfaces with three principal curvatures, {\em Math. Z.} {\bf 187} (1984),
433--452.

\bibitem{Mi8} Miyaoka R., Taut embeddings and Dupin hypersurfaces, in Dif\/ferential Geometry of Submanifolds (Kyoto, 1984),
Editor K. Kenmotsu, {\it Lecture Notes in Math.}, Vol.~1090,
Springer, Berlin~-- New York, 1984, 15--23.

\bibitem{Mi2} Miyaoka R., Dupin hypersurfaces and a Lie invariant, {\em Kodai Math. J.} {\bf 12} (1989), 228--256.

\bibitem{Mi3} Miyaoka R., Dupin hypersurfaces with six principal curvatures, {\em Kodai Math. J.} {\bf 12} (1989),
308--315.

\bibitem{Mi6} Miyaoka R., The linear isotropy group of $G_2/SO(4)$, the Hopf f\/ibering and isoparametric
hypersurfaces, {\em Osaka J. Math.} {\bf 30} (1993), 179--202.

\bibitem{Mi10} Miyaoka R., The Dorfmeister--Neher's theorem on isoparametric hypersurfaces,
\href{http://arxiv.org/abs/math.DG/0602519}{math.DG/0602519}.

\bibitem{MO} Miyaoka  R., Ozawa T., Construction of taut
embeddings and Cecil-Ryan conjecture, in Geometry of Manifolds, Editor K. Shiohama, {\it Perspect. Math.}, Vol.~8,
Academic Press, Boston, 1989, 181--189.

\bibitem{MC} Morse M., Cairns S., Critical point theory
in global analysis and dif\/ferential topology, Academic Press, New York, 1969.

\bibitem{Mul} Mullen S., Isoparametric systems on symmetric spaces, in Geometry and Topology of Submanifolds~VI,
Editors F. Dillen et al.,
World Scientif\/ic, River Edge, NJ, 1994, 152--154.

\bibitem{Mu} M\"{u}nzner H.-F., Isoparametrische Hyperf\/l\"{a}chen in Sph\"{a}ren, {\em Math. Ann.} {\bf 251} (1980), 57--71.

\bibitem{Mu2} M\"{u}nzner H.-F., Isoparametrische Hyperf\/l\"{a}chen in Sph\"{a}ren II: \"{U}ber die Zerlegung
der Sph\"{a}re in Ballb\"{u}ndel, {\em Math. Ann.} {\bf 256} (1981), 215--232.

\bibitem{N1} Niebergall R., Dupin hypersurfaces in ${\bf R}^5$.~I, {\em Geom. Dedicata} {\bf 40} (1991), 1--22.

\bibitem{N2} Niebergall R., Dupin hypersurfaces in ${\bf R}^5$.~II, {\em Geom. Dedicata} {\bf 41} (1992), 5--38.

\bibitem{NieR1} Niebergall  R., Ryan P., Isoparametric hypersurfaces -- the af\/f\/ine case, in Geometry
and Topology of Submanifolds V, Editors F. Dillen et al., World Scientif\/ic,
River Edge, NJ, 1993, 201--214.

\bibitem{NieR2} Niebergall  R., Ryan P., Af\/f\/ine isoparametric hypersurfaces, {\em Math. Z.} {\bf 217} (1994),
479--485.

\bibitem{NieR3} Niebergall  R., Ryan P., Focal sets in af\/f\/ine geometry, in  Geometry and Topology of
Submanifolds VI, Editors F.~Dillen et al., World Scientif\/ic,
River Edge, NJ, 1994, 155--164.

\bibitem{NieR4} Niebergall  R., Ryan P., Af\/f\/ine Dupin surfaces, {\em Trans. Amer. Math. Soc.} {\bf 348} (1996),
1093--1117.

\bibitem{NieR} Niebergall  R., Ryan P., Real hypersurfaces in complex space forms,
in Tight and Taut Submanifolds (Berkeley, CA, 1994), Editors T.~Cecil and
S.-S.~Chern, {\it Math. Sci. Res. Inst. Publ.}, Vol.~32, Cambridge Univ. Press, Cambridge, 1997, 233--305.

\bibitem{Nom1} Nomizu K., Some results in E. Cartan's
theory of isoparametric families of hypersurfaces, {\em Bull. Amer. Math. Soc.} {\bf 79} (1973), 1184--1188.

\bibitem{Nom2} Nomizu K., \'{E}lie Cartan's work on isoparametric families of hypersurfaces,
in Dif\/ferential geometry (Stanford, 1973), Editors S.-S.~Chern and R.~Osserman, {\it Proc.
Sympos. Pure Math.}, Vol.~27, Amer. Math. Soc.,
Providence, RI, 1975, Part 1, 191--200.

\bibitem{Nom3} Nomizu K., On isoparametric hypersurfaces in Lorentzian space forms,
{\em Japan. J. Math.} {\bf 7} (1981), 217--226.

\bibitem{Olm} Olmos C., Isoparametric submanifolds and their
homogeneous structure, {\em J. Differential Geom.} {\bf 38} (1993), 225--234.


\bibitem{Oz} Ozawa T., On the critical sets of distance functions
to a taut submanifold, {\em Math. Ann.} {\bf 276} (1986), 91--96.

\bibitem{OT} Ozeki  H., Takeuchi M., On some types of
isoparametric hypersurfaces in spheres.~I, {\em T\^{o}hoku Math. J.} {\bf 27} (1975), 515--559.

\bibitem{OT2} Ozeki  H., Takeuchi M., On some types of isoparametric hypersurfaces in spheres.~II, {\em T\^{o}hoku
Math. J.} {\bf 28} (1976), 7--55.

\bibitem{PalT1} Palais  R., Terng C.-L., A general theory
of canonical forms, {\em Trans. Amer. Math. Soc.} {\bf 300} (1987), 771--789.

\bibitem{PalT} Palais  R., Terng C.-L., Critical point
theory and submanifold geometry, {\it Lecture Notes in Math.}, Vol.~1353,
Springer, Berlin~-- New York, 1988.

\bibitem{Peng} Peng C.-K., Hou Z., A remark on the
isoparametric polynomials of degree~6,
in Dif\/ferential Geometry and Topology (Tianjin, 1986--87), Editors B.~Jiang et al.,
{\it Lecture Notes in Math.}, Vol~1369, Springer, Berlin~-- New York, 1989, 222--224.

\bibitem{P1} Pinkall U., Dupin'sche Hyperf\/l\"{a}chen,
Dissertation, Univ. Freiburg, 1981.

\bibitem{P6} Pinkall U., Letter to T. Cecil, December~5, 1984.

\bibitem{P2} Pinkall U., Dupin'sche Hyperf\/l\"{a}chen in $E^4$, {\em Manuscripta Math.} {\bf 51} (1985), 89--119.

\bibitem{P3} Pinkall U., Dupin hypersurfaces, {\em Math. Ann.} {\bf 270} (1985), 427--440.

\bibitem{P4} Pinkall U., Curvature properties of taut submanifolds, {\em Geom. Dedicata} {\bf 20} (1986), 79--83.

\bibitem{PT1} Pinkall U., Thorbergsson G., Deformations of Dupin hypersurfaces,
{\em Proc. Amer. Math. Soc.} {\bf 10} (1989), 1037--1043.

\bibitem{PT3} Pinkall  U., Thorbergsson  G., Examples of
inf\/inite dimensional isoparametric submanifolds, {\em Math. Z.} {\bf 205} (1990), 279--286.

\bibitem{Pr1} Pratt M.J., Cyclides in computer aided geometric design, {\em Comput.
Aided Geom. Design} {\bf 7} (1990), 221--242.

\bibitem{Pr2} Pratt M.J., Cyclides in computer aided geometric design.~II, {\em Comput.
Aided Geom. Design} {\bf 12} (1995), 131--152.

\bibitem{Reck} Reckziegel H., On the eigenvalues of the shape
operator of an isometric immersion into a space of constant curvature,
{\em Math. Ann.} {\bf 243} (1979), 71--82.

\bibitem{RRT} Riveros C.M.C., Rodrigues L.A., Tenenblat K., On Dupin hypersurfaces with constant M\"{o}bius
curvature, {\em Pacific J. Math.} {\bf 236} (2008), 89--103.

\bibitem{RT} Riveros  C.M.C., Tenenblat K., On four dimensional Dupin hypersurfaces in Euclidean space,
{\em An. Acad. Brasil Ci\^{e}nc.} {\bf 75} (2003), 1--7.

\bibitem{RT2} Riveros  C.M.C., Tenenblat K., Dupin hypersurfaces in ${\bf R}^5$,
{\em Canad. J. Math.} {\bf 57} (2005), 1291--1313.

\bibitem{RoT} Rodrigues  L.A., Tenenblat K., A characterization of Moebius isoparametric hypersurfaces of the sphere,
Preprint, 2008.

\bibitem{Ryan3} Ryan P.J., Hypersurfaces with parallel Ricci tensor, {\em Osaka J. Math.} {\bf 8} (1971), 251--259.

\bibitem{Schrott} Schrott  M., Odehnal B., Ortho-circles of Dupin cyclides, {\em J. Geom. Graph.} {\bf 10} (2006), 73--98.

\bibitem {Seg1} Segre B., Una propriet\'{a} caratteristica de tre sistemi $\infty^1$ di superf\/icie, {\em Atti Acc. Sc. Torino}
{\bf 59} (1924), 666--671.

\bibitem {Seg} Segre B., Famiglie di ipersuperf\/icie
isoparametrische negli spazi euclidei ad un qualunque numero
di demesioni, {\em Atti. Accad. naz Lincie Rend. Cl. Sci. Fis.
Mat. Natur.} {\bf 27} (1938), 203--207.

\bibitem{Sin} Singley D., Smoothness theorems for the principal
curvatures and principal vectors of a hypersurface, {\em Rocky
Mountain J. Math.} {\bf 5} (1975), 135--144.

\bibitem{Sol1} Solomon B., The harmonic analysis of cubic isoparametric minimal hypersurfaces.~I. Dimensions 3 and 6,
{\em Amer. J. Math.} {\bf 112} (1990), 157--203.

\bibitem{Sol2} Solomon B., The harmonic analysis of cubic isoparametric minimal hypersurfaces.~II. Dimensions 12 and~24,
{\em Amer. J. Math.} {\bf 112} (1990), 205--241.

\bibitem{Sol3} Solomon B., Quartic isoparametric hypersurfaces and quadratic forms, {\em Math. Ann.}
{\bf 293} (1992), 387--398.

\bibitem{Som} Somigliana C., Sulle relazione fra il principio di Huygens e l'ottica geometrica, {\em Atti Acc. Sc. Torino}
{\bf 54} (1918--1919), 974--979 (see also in {\em Memorie Scelte}, 434--439).

\bibitem{Sp} Spivak M., A comprehensive introduction to dif\/ferential geometry, Vol.~4, Publish or Perish,
Boston, 1975.

\bibitem{SD1} Srinivas  Y.L., Dutta D., Blending and joining using cyclides, {\em ASME Trans. J.
Mechanical Design} {\bf 116} (1994), 1034--1041.

\bibitem{SD2} Srinivas  Y.L., Dutta D.,
An intuitive procedure for constructing complex objects using
cyclides, {\em Computer-Aided Design} {\bf 26} (1994), 327--335.

\bibitem{SD3} Srinivas  Y.L., Dutta D.,
Cyclides in geometric modeling: computational tools for
an algorithmic infrastructure, {\em ASME Trans. J. Mechanical
Design} {\bf 117} (1995), 363--373.

\bibitem{SD4} Srinivas  Y.L., Dutta D.,
Rational parametrization of parabolic cyclides, {\em Comput.
Aided Geom. Design} {\bf 12} (1995), 551--566.

\bibitem{Stolz} Stolz S., Multiplicities of Dupin hypersurfaces,
{\em Invent. Math.} {\bf 138} (1999), 253--279.

\bibitem{Str} Str\"{u}bing W., Isoparametric submanifolds,
{\em Geom. Dedicata} {\bf 20} (1986), 367--387.

\bibitem{Takagi} Takagi R., A class of hypersurfaces with constant principal
curvatures in a sphere, {\em J. Differential Geom.} {\bf 11} (1976), 225--233.

\bibitem{TT} Takagi R., Takahashi T., On the principal
curvatures of homogeneous hypersurfaces in a sphere, in Dif\/ferential
Geometry (in Honor of Kentaro Yano), Kinokuniya, Tokyo, 1972, 469--481.

\bibitem{Tak} Takeuchi M., Proper Dupin hypersurfaces
generated by symmetric submanifolds, {\em Osaka Math. J.} {\bf 28} (1991), 153--161.

\bibitem{TK} Takeuchi  M., Kobayashi S., Minimal imbeddings
of $R$-spaces, {\em J. Differential Geom.} {\bf 2} (1968), 203--215.

\bibitem{Tang} Tang Z.-Z., Isoparametric hypersurfaces
with four distinct principal curvatures, {\em Chinese Sci. Bull.} {\bf 36}
(1991), 1237--1240.

\bibitem{Tang2} Tang Z.-Z., Multiplicities of equifocal hypersurfaces in symmetric spaces, {\em Asian J. Math.} {\bf 2} (1998), 181--214.

\bibitem{Te1} Terng C.-L., Isoparametric submanifolds and their
Coxeter groups, {\em J. Differential Geom.} {\bf 21} (1985), 79--107.

\bibitem{Te4} Terng C.-L., Convexity theorem
for isoparametric submanifolds, {\em Invent. Math.} {\bf 85} (1986), 487--492.

\bibitem{Te3} Terng C.-L., Submanifolds with
f\/lat normal bundle, {\em Math. Ann.} {\bf 277} (1987), 95--111.

\bibitem{Te5} Terng C.-L., Proper Fredholm submanifolds of Hilbert space,
{\em J. Differential  Geom.} {\bf 29} (1989), 9--47.

\bibitem{Te2} Terng C.-L., Recent progress in
submanifold geometry, in Dif\/ferential Geometry: Partial Dif\/ferential Equations on Manifolds (Los Angeles, 1990),
Editors R.~Greene and S.T.~Yau,
{\it Proc. Sympos. Pure Math.}, Vol.~54,
Amer. Math. Soc., Providence, RI, 1993, Part 1, 439--484.

\bibitem{TTh} Terng C.-L., Thorbergsson G., Submanifold
geometry in symmetric spaces, {\em J. Differential Geom.} {\bf 42} (1995), 665--718.

\bibitem{TTh1} Terng  C.-L., Thorbergsson G., Taut immersions
into complete Riemannian manifolds, in Tight and Taut Submanifolds, Editors T.~Cecil and
S.-S.~Chern, {\it Math. Sci. Res. Inst. Publ.}, Vol.~32, Cambridge Univ. Press, Cambridge, 1997, 181--228.

\bibitem{Th1} Thorbergsson G., Dupin hypersurfaces, {\em Bull. London Math. Soc.} {\bf 15} (1983), 493--498.

\bibitem{Th2} Thorbergsson G., Isoparametric
foliations and their buildings, {\em Ann. Math.} {\bf 133} (1991), 429--446.

\bibitem{Th6} Thorbergsson G., A survey
on isoparametric hypersurfaces and their generalizations,
in  Handbook of Dif\/ferential Geometry, Vol.~I, Editors F.~Dillen and L.~Verstraelen,
Elsevier Science, Amsterdam, 2000, 963--995.

\bibitem{Toben} T\"{o}ben D., Parallel focal structure and singular Riemannian foliations,
{\em Trans. Amer. Math. Soc.} {\bf 358} (2006), 1677--1704, \href{http://arxiv.org/abs/math.DG/0403050}{math.DG/0403050}.

\bibitem{Ver} Verh\'{o}czki L., Isoparametric submanifolds
of general Riemannian manifolds, in Dif\/ferential Geometry and
Its Applications (Eger, 1989), Editors J.~Szenthe and L.~Tam\'{a}ssy, {\it Colloq. Math. Soc.
J\'{a}nos Bolyai}, Vol.~56, North-Holland, Amsterdam, 1992, 691--705.

\bibitem{Wc} Wang C.-P., Surfaces in M\"{o}bius geometry,
{\em Nagoya Math. J.} {\bf 125} (1992), 53--72.

\bibitem{Wc1} Wang C.-P., M\"{o}bius geometry for hypersurfaces in $S^4$, {\em Nagoya Math. J.} {\bf 139} (1995),
1--20.

\bibitem{Wc2} Wang C.-P., M\"{o}bius geometry of submanifolds in $S^n$, {\em Manuscripta Math.} {\bf 96} (1998),
517--534.

\bibitem{Wang3} Wang Q.-M., Isoparametric maps of Riemannian manifolds and their applications, in Advances
in Science of China, Mathematics, Vol.~2, Editors C.H.~Gu and Y. Wang, Wiley-Interscience,
New York, 1986, 79--103.

\bibitem{Wang1} Wang Q.-M., Isoparametric functions on
Riemannian manifolds.~I, {\em Math. Ann.} {\bf 277} (1987), 639--646.

\bibitem{Wang2} Wang Q.-M., On the topology
of Clif\/ford isoparametric hypersurfaces, {\em J. Differential Geom.} {\bf 27} (1988), 55--66.

\bibitem{West} West A., Isoparametric systems, in Geometry and
Topology of Submanifolds, Editors J.-M.~Morvan and L.~Verstraelen, World
Scientif\/ic, River Edge, NJ, 1989, 222--230.

\bibitem{West1} West A., Isoparametric systems on symmetric spaces, in Geometry and Topology of
Submanifolds V, Editors F.~Dillen et al., World Scientif\/ic, River Edge, NJ, 1993, 281--287.

\bibitem{Wu} Wu B., Isoparametric submanifolds of hyperbolic
spaces, {\em Trans. Amer. Math. Soc.} {\bf 331} (1992), 609--626.

\bibitem{Wu2} Wu B., A f\/initeness theorem for isoparametric hypersurfaces, {\em Geom. Dedicata} {\bf 50} (1994),
247--250.

\bibitem{Yau} Yau S.-T., Open problems in geometry,
in Dif\/ferential Geometry: Partial Dif\/ferential Equations on Manifolds (Los Angeles, 1990),
Editors R.~Greene and S.T.~Yau,
{\it Proc. Sympos. Pure Math.}, Vol.~54,
Amer. Math. Soc., Providence, RI, 1993, Part 1, 439--484.

\bibitem{Zhao} Zhao Q., Isoparametric submanifolds of hyperbolic
spaces, {\em Chinese J. Contemp. Math.} {\bf 14} (1993), 339--346.

\end{thebibliography}
\end{document}